\documentclass[12pt]{amsart}
\usepackage{latexsym,amsmath,amssymb,amscd,epsfig,color}
\usepackage[pdftex]{hyperref} %for Pdf Latex

\newtheorem{thm}{Theorem}[section]
\newtheorem{lemma}[thm]{Lemma}
\newtheorem{proposition}[thm]{Proposition}
\newtheorem{coroll}[thm]{Corollary}
\newtheorem{corollary}[thm]{Corollary}
\theoremstyle{definition}
\newtheorem{definition}[thm]{Definition}

\theoremstyle{remark}
\newtheorem{remark}[thm]{Remark}

\newtheorem{example}[thm]{Example}
\newtheorem{examples}[thm]{Examples}

\newcommand{\Z}{{\mathbb Z}}

\newcommand{\SL}{{\rm SL}}

\newcommand{\xac}{{\mathcal{X}_{ac}}}
\newcommand{\mac}{{\mathcal{M}_{ac}}}
\newcommand{\kh}{{\scriptstyle{\mathbf H}}}
\newcommand{\cF}{{\mathfrak{F}}}
\newcommand{\khF}{{\kh\cF}}
\newcommand{\cC}{{\mathfrak{C}}}
\newcommand{\khC}{{\kh\cC}}

\newcommand{\cA}{\mathcal{A}}

\newcommand{\Hl}{\{ H_\lambda \}_{\lambda \in \Lambda}}

\newcommand{\bigast}{\ensuremath{\mathop{{\ast}}}}
\def\bfg{{\bf G}}
\def\cals{{\mathcal S}}
\def\calsplus{{{\mathcal S}^+}}
\def\mapdown#1{\Big\downarrow\rlap{$\vcenter{\hbox{$\scriptstyle#1$}}$}}
\def\mapright#1{\smash{\mathop{\longrightarrow}\limits^{#1}}}
\def\min{{\rm min}}
\def\pra{{\par}}
\def\span{{\rm span}}
\def\tilG{{\widetilde G}}

\newcommand{\N}{\mathbb{N}}

\title{Infinite groups with fixed point properties}

\author[Arzhantseva]{G. Arzhantseva}
\address[G. Arzhantseva]{Universit\'{e} de Gen\`{e}ve,
Section de Math\'{e}matiques, 2-4 rue du Li\`{e}vre, Case postale
64, 1211 Gen\`{e}ve 4, Switzerland}
\email{Goulnara.Arjantseva@math.unige.ch}

\author[Bridson]{M.R. Bridson}
\address[M.R. Bridson]{Mathematical Institute, 24-29 St Giles', Oxford, UK}
\email{bridson@maths.ox.ac.uk}

\author[Januszkiewicz]{T. Januszkiewicz}
\address[T. Januszkiewicz]{
Department of Mathematics,
The Ohio State University, 231 W. 18th Ave., Columbus, OH 43210,
USA and The Mathematical Institute of Polish Academy of
Sciences. On leave from Instytut Matematyczny, Uniwersytet Wroc\l
awski} \email{tjan@math.ohio-state.edu}

\author[Leary]{I.J. Leary}
\address[I.J. Leary]{Department of Mathematics,
The Ohio State University, 231 W. 18th Ave., Columbus, OH 43210,
USA} \email{leary@math.ohio-state.edu}

\author[Minasyan]{A. Minasyan}
\address[A. Minasyan]{%Universit\'{e} de Gen\`{e}ve,
School of Mathematics,
University of Southampton, Southampton, SO17 1BJ, United
Kingdom.}  \email{aminasyan@gmail.com}

\author[\'Swi\c atkowski]{J.~\'Swi\c atkowski}
\address[J. \'Swi\c atkowski]{Instytut Matematyczny, Uniwersytet
  Wroclawski, pl. Grunwaldzki 2/4, 50-384 Wroclaw, Poland}
\email{swiatkow@math.uni.wroc.pl}

\thanks{The work of Arzhantseva and Minasyan and a visit to Geneva
  by \'Swi\c atkowski were supported by Swiss National Science
  Foundation (SNSF) grants PP002-68627, PP002-116899/1; visits of Bridson
  to Geneva were also supported by the SNSF. The work of
  \'Swi\c atkowski was supported by the Polish Ministry of Science,
  Higher Education (MNiSW) grant N201 012 32/0718, and a visit by
  him to Ohio was supported by the Ohio State University's
  Mathematics Research Institute.
  Bridson was supported by a Senior Fellowship from the EPSRC of the
  United Kingdom and a Royal
  Society Wolfson Research Merit Award.  Januszkiewicz was partially
  supported by NSF grant DMS-0706259.  Leary was partially supported
  by NSF grant DMS-0505471 and by the Heilbronn Institute.
  Several of the authors worked on this
  manuscript at MSRI, Berkeley, where research
  is supported in part by NSF grant DMS-0441170.}

\dedicatory{Dedicated to  Michael W. Davis on the occasion of his 60th birthday.}
\begin{document}
\begin{abstract} We construct finitely generated groups with
strong fixed point properties.
Let $\xac$ be the class of Hausdorff spaces of finite
  covering dimension which are mod-$p$ acyclic for at least one prime
  $p$.  We produce the first examples of infinite
  finitely generated groups $Q$ with the property that for
  any action of $Q$ on any $X\in \xac$, there is a global fixed
  point.  Moreover, $Q$ may be chosen to be simple
  and to have Kazhdan's property (T).  We construct a finitely
  presented infinite group $P$ that admits no non-trivial action
  on any manifold in $\xac$.  In building
  $Q$, we exhibit new families of hyperbolic groups:
  for each $n\geq
  1$ and each prime $p$, we construct a non-elementary hyperbolic
  group $G_{n,p}$ which has a generating set of size $n+2$, any proper
  subset of which generates a finite $p$-group.
\end{abstract}
\subjclass[2000]{Primary 20F65, 20F67 Secondary  57S30, 55M20}
\keywords{Acyclic spaces, Kazhdan's property (T), relatively hyperbolic groups, simplices of groups}

\maketitle
%%%%%%%%%%%%%%%%%%%%%%%%%%%%%%%%%%%%%%%%%%%%%%%%%%%%%
\section{Introduction}
%%%%%%%%%%%%%%%%%%%%%%%%%%%%%%%%%%%%%%%%%%%%%%%%%%%%%

We present three templates for proving fixed point theorems; two
are based on relative small cancellation theory and one is based
on the Higman Embedding Theorem. Each template demands as input a
sequence of groups with increasingly strong fixed point
properties. By constructing such sequences we prove the following
fixed point theorems.
%Before stating the theorems, we define some classes of spaces.

%\begin{definition}\label{def:xac}
For $p$ a prime, one says that a space is \emph{mod-$p$ acyclic} if it has the
same mod-$p$ \v Cech cohomology as a point.  Let $\xac$ be the class
of all Hausdorff spaces $X$ of finite covering dimension such that
there is a prime $p$ for which $X$ is mod-$p$ acyclic.  Let $\mac$
denote the subclass of manifolds in $\xac$. %\end{definition}

Note that the class $\xac$ contains all finite dimensional
contractible spaces and all finite dimensional
contractible CW-complexes.

%Recall that a group $Q$ is called {\it
%periodic} if each element of $Q$ has a finite order.

\begin{thm}\label{thm:fp}
There is an infinite
finitely generated group $Q$ that cannot act without a global fixed point
on any $X\in \mathcal{X}_{ac}$.  If $X\in \xac$ is mod-$p$ acyclic,
then so is the fixed point set for any action of $Q$ on $X$.
For any countable group $C$, the group $Q$ can be chosen to have either
the additional properties (i), (ii) and (iii) or (i), (ii) and
(iii$)'$ described below:
\begin{itemize}
  \item[(i)] $Q$ is simple;
  \item[(ii)] $Q$ has Kazhdan's property (T);
  \item[(iii)] $Q$ contains an isomorphic copy of $C$;
  \item[(iii$)'$] $Q$ is periodic.
\end{itemize}
\end{thm}

Since a countable group can contain only countably many finitely
generated subgroups, it follows from property~(iii) that there are
continuously many (i.e., $2^{\aleph_0}$) non-isomorphic groups $Q$
with the fixed point property described in Theorem~\ref{thm:fp}.

Recall that for countable groups, Kazhdan's property (T) is equivalent to the
fact that every isometric action of the group on a Hilbert space
has a global fixed point.

No non-trivial finite group has a fixed point property as strong as the one in
Theorem~\ref{thm:fp}.
Any finite group not of prime power order acts without a global fixed
point on some finite dimensional contractible simplicial complex.
%If the assertion above needs a reference, IJL's paper listed as
%\cite{L} will do.  If the assertion below (about $p$-groups) needs a
%reference, then Bredon's book listed as \cite{bredon} will do.
Smith theory tells us that the fixed point set for any action of a
finite $p$-group on a finite dimensional mod-$p$ acyclic space is
itself mod-$p$ acyclic, but it is easy to construct an action of a
non-trivial finite $p$-group on a 2-dimensional mod-$q$ acyclic space
without a global fixed point if $q$ is any prime other than $p$.  Since
the fixed point property of Theorem~\ref{thm:fp} passes to quotients,
it follows that none of the groups $Q$ can admit a non-trivial finite
quotient. This further restricts the ways in which $Q$ can act on acyclic
spaces. For example, if $X\in \xac$ is a locally finite simplicial
complex and $Q$ is acting simplicially, then the action of $Q$
on the successive star neighbourhoods ${\rm{st}}_{n+1}(x) := {\rm{st}}({\rm{st}}_{n}(x))$
of a fixed point $x\in X$ must be trivial, because $ {\rm{st}}_n(x)$ is $Q$-invariant
and there is no non-trivial map from $Q$ to the finite group ${\rm Aut}( {\rm{st}}_n(x))$.
Since $X= \bigcup_n  {\rm{st}}_n(x)$, we deduce:

\begin{corollary} \label{cor:simp} None of the groups
$Q$ from Theorem~\ref{thm:fp} admit a non-trivial simplicial action on
any locally-finite simplicial complex $X\in \xac$. (We do not assume that
$Q$ satisfies any of the conditions (i) to (iii)$'$.)
\end{corollary}

The ideas behind the above corollary can be further developed in different
directions. The following was suggested to us by N. Monod:

\begin{corollary} \label{cor:} None of the groups
$Q$ from Theorem~\ref{thm:fp} admit
a non-trivial isometric action on
any proper metric space $X\in \xac$. (Again, we do not assume that
$Q$ satisfies any of the conditions (i) to (iii)$'$.)
\end{corollary}

Indeed, suppose that $x_0\in X$ is fixed by the action of $Q$. Then we
have a homomorphism $\phi:Q\to {\rm Isom}(X)$ such that $\phi(Q) \le
{\rm{Stab}}(x_0)$, where ${\rm{Stab}}(x_0)$
 denotes the stabilizer of $x_0$ in the full isometry group
  of $X$.  Since $X$ is proper, by a
general version of Arzel\`a-Ascoli theorem (see \cite[I.3.10]{BH}) the
group ${\rm{Stab}}(x_0)$ is compact (when endowed with the topology of uniform
convergence on compact subsets). If the action of $Q$ on $X$
were non-trivial, then $\phi(g)\in {\rm{Stab}}(x_0)$ would be non-trivial
 for some $g \in Q$, and by the Peter-Weyl theorem (see \cite[Sec. 4]{Rob})
  there would be  a
finite-dimensional irreducible representation $\rho$ of ${\rm{Stab}}(x_0)$ such
that $\rho(\phi(g))\neq {\rm Id}$. But then $\rho \circ \phi(Q)$ would be a
non-trivial finitely generated linear group, which is a contradiction
because such groups are residually
finite, by a well-known result of Malcev (see
\cite[III.7.11]{L-S}), whereas  $Q$ has no non-trivial
finite quotients.

Using similar arguments we can also rule out non-trivial real analytic
actions of the groups $Q$ (from Theorem~\ref{thm:fp}) on any acyclic
manifold $M$: since the action fixes a point $x_0 \in M$, the image of
every element of $Q$ in the group ${\rm Homeo}(M)$ can be regarded as
an $n$-tuple of formal power series on $n$ variables with real
coefficients and trivial constant terms, where $n = \dim(M)$. Thus,
such an action gives rise to a homomorphism from $Q$ into the group
$H_n$ of such invertible $n$-tuples of formal power series with
respect to composition. One can argue that $H_n$ is embedded into an
inverse limit of linear groups, therefore every finitely generated
subgroup of $H_n$ is residually finite.  As before, the absence of
non-trivial finite quotients for $Q$ yields the triviality of the
above homomorphism.

A stronger result concerning triviality of
actions on manifolds can be obtained more directly from Theorem~\ref{thm:fp}
by using properties (i) and (iii).

\begin{proposition} \label{prop:simp}
  A simple group $G$ that contains, for each $n>0$ and each prime $p$, a
  copy of $(\Z_p)^n$ admits no non-trivial action
  on any $X\in \mac$.  The group $Q$ in Theorem~\ref{thm:fp} may be
  chosen to have this property.
\end{proposition}

Finite $p$-groups have global fixed points whenever they act on
compact Hausdorff spaces that are mod-$p$ acyclic, but the groups $Q$ do not
have this property. Indeed, if $Q$ is infinite and has property (T) then
it will be non-amenable, hence the natural action of $Q$ on
the space of finitely-additive
probability measures on $Q$ will not have a global fixed point, and this
space is compact, contractible, and Hausdorff.

We know of no finitely presented group enjoying the fixed point
property described in Theorem~\ref{thm:fp}. However, using techniques quite different from those used to construct
the groups $Q$, we shall exhibit finitely presented groups that cannot
act on a range of spaces. In particular we construct groups of the following
type.

\begin{thm} \label{thm:triv} There exist finitely presented infinite
groups $P$ that have no non-trivial action on any manifold $X\in \mac$.
\end{thm}

Certain of the Higman-Thompson groups \cite{higman}, such as
R.~Thompson's vagabond group $V$, can also serve in the role
of $P$ (cf.~Remark~\ref{rem:last}).

%%%%%%%%%%%%%%%%%% Ian's intro %%%%%%%%%%%%%%%%%%%%%%%%%%%%%%
Theorem~\ref{thm:fp}
answers a question of P.~H.~Kropholler, who asked whether
there exists a countably infinite group $G$ for which every
finite-dimensional contractible $G$-CW-complex has a global
fixed point.
This question is motivated by Kropholler's study of the closure
operator $\kh$ for classes of groups, and by the class $\khF$ obtained
by applying this operator to the class $\cF$ of all finite groups
\cite{krop}.  Briefly, if $\cC$ is a class of groups, then the class
$\khC$ is the smallest class of groups that contains $\cC$ and has the
property that if the group $G$ admits a
finite-dimensional contractible $G$-CW-complex $X$ with all stabilizers
already in $\khC$, then $G$ is itself in $\khC$.  Kropholler showed
that any torsion-free group in $\khF$ of type $FP_\infty$ has finite
cohomological dimension.  Since Thompson's group $F$ is torsion-free
and of type $FP_\infty$ but has infinite cohomological
dimension~\cite{brogeo}, it follows that $F$ is not in $\khF$.  Until
now, the only way known to show that a group is not in $\khF$ has been
to show that it contains Thompson's group as a subgroup.  If $Q$ is
any of the groups constructed in Theorem~\ref{thm:fp}, then $Q$ has
the property that for any class $\cC$ of groups, $Q$ is in the class
$\khC$ if and only if $Q$ is already in the class $\cC$.  In
particular, $Q$ is not in the class $\khF$.  Note that many of the
groups constructed in Theorem~\ref{thm:fp} cannot contain Thompson's
group $F$ as a subgroup, for example the periodic groups.

Our strategy for proving Theorems \ref{thm:fp}~and~\ref{thm:triv} is
very general.  First, we express our class of spaces as a countable
union $\mathcal{X}=\cup_{n\in \N} \mathcal{X}_n$.  For instance,
if all spaces in $\mathcal{X}$ are finite-dimensional, then
$\mathcal{X}_n$ may be taken to consist of all $n$-dimensional spaces
in $\mathcal{X}$.
Secondly, we construct finitely generated groups $G_n$ that have the
required properties for actions on any $X\in\mathcal{X}_n$.
Finally, we apply the templates described below to produce the required
groups.

\medskip
\noindent\emph{Template {\bf FP}: ruling out fixed-point-free actions.} If
there is a sequence of finitely generated \emph{non-elementary
relatively hyperbolic} groups $G_n$ such that $G_n$ cannot act
without a fixed point on any $X\in\mathcal{X}_n$, then there is an infinite
finitely generated group that cannot act without a fixed-point on
any $X\in \mathcal{X}$.

\medskip

\noindent\emph{Template $\text{\bf NA}_{fg}$: ruling out non-trivial
  actions.} If there is a sequence of non-trivial finitely generated
  groups $G_n$ such that $G_n$ cannot act non-trivially on any $X\in
  \mathcal{X}_n$, then there is an infinite finitely generated group
  that cannot act non-trivially on any $X\in \mathcal{X}$.

\medskip

\noindent\emph{Template $\text{\bf NA}_{fp}$: finitely presented
groups that cannot act.} Let $(G_n;\xi_{n,j})\ (n\in\N,
j=1,\dots,J)$ be a recursive system of non-trivial groups and monomorphisms
$\xi_{n,j}:G_n\to G_{n+1}$.
Suppose that each $G_{n+1}$ is generated by
$\bigcup_j\xi_{n,j}(G_n) $ and that for every $m\in\N$ there
exists $n\in\N$ such that
 $G_n$ cannot act non-trivially on any
$X\in \mathcal{X}_m$.  Then there exists an infinite
finitely presented group that cannot act non-trivially on any $X\in
\mathcal{X}$.

\medskip
Only the first and third templates are used in the construction of the
groups $P$~and~$Q$.  We include the second template with a view to
further applications.

The engine that drives the first two templates is the existence of
common quotients established in Theorem~\ref{thm:comm_quot} below.
The proof of this theorem, given in Section~\ref{sec:qq}, is based on
the following result of Arzhantseva, Minasyan and Osin~\cite{SQ},
obtained using small cancellation theory over relatively hyperbolic
groups: any two finitely generated non-elementary relatively
hyperbolic groups $G_1, G_2$ have a common non-elementary relatively
hyperbolic quotient $H$.

\begin{thm} \label{thm:comm_quot}
Let $\{G_n\}_{n \in \N}$ be a countable collection of finitely
generated non-ele\-men\-tary relatively hyperbolic groups. Then
 there exists an infinite finitely generated group $Q$
that is a quotient of $G_n$ for every $n \in \N$.

Moreover,  if $C$ is an arbitrary countable group, then such a group $Q$
can be made to satisfy the following conditions
\begin{itemize}
\item[(i)] $Q$ is a simple group;
\item[(ii)] $Q$ has Kazhdan's
property (T);
\item[(iii)] $Q$ contains an isomorphic copy of $C$.
\end{itemize}

If the $G_n$ are non-elementary word hyperbolic groups, then
claim {\rm (iii)} above can be replaced by
\begin{itemize}
 \item[(iii$)'$] $Q$ is periodic.
\end{itemize}
\end{thm}

This result immediately implies the validity of templates \textbf{FP}
and $\textbf{NA}_{fg}$.  Indeed, if $G_n$ are the hypothesized groups
of template \textbf{FP}, the preceding theorem furnishes us with a
group $Q$ that, for each $n\in\mathbb N$, is a quotient of $G_n$ and
hence cannot act without a fixed point on any $X\in \mathcal{X}_n$.
Now let $G_n$ be the hypothesized groups of template
$\textbf{NA}_{fg}$.  They are not assumed to be relatively hyperbolic.
We consider groups $A_n:=G_n\ast G_n \ast G_n$, which also cannot act
non-trivially on any $X\in \mathcal{X}_n$.  The group $A_n$ is
non-elementary and relatively hyperbolic as a free product of three
non-trivial groups.  Therefore, Theorem~\ref{thm:comm_quot} can be
applied to the sequence of groups $A_n$, providing a group $Q_1$
which, as a quotient of $A_n$, cannot act non-trivially on any
$X\in\mathcal X_n$ for any $n\in\mathbb N$.

Following the above strategy to prove Theorem~\ref{thm:fp}, we
first represent $\mathcal{X}_{ac}$ as a countable union
$\mathcal{X}_{ac}=\cup_{n,p}\mathcal{X}_{n,p}$, where, for each
prime number $p$, the class $\mathcal{X}_{n,p}$ consists of all
mod-$p$ acyclic spaces of dimension $n$. Then, in
Section~\ref{sec:retract}, we construct the groups required by
template \textbf{FP}, proving the following result.

\begin{thm} \label{thm:retract_gp}
For each $n\in \N$ and every prime $p$, there exists a non-elementary
word hyperbolic group $G_{n,p}$ such that any action of $G_{n,p}$ by
homeomorphisms on any space $X\in \mathcal{X}_{n,p}$ has the property
that the global fixed point set is mod-$p$ acyclic (and in particular
non-empty).
\end{thm}

The mod-$p$ acyclicity of the fixed point for the action of $G_{n,p}$
on the space $X$ is a consequence of the following {\it
$(n,p)$-generation property}: there is a generating set $S$ of
$G_{n,p}$ of cardinality $n+2$ such that any proper subset of $S$
generates a finite $p$-subgroup.

For certain small values of the parameters examples of non-elementary
word hyperbolic groups with the $(n,p)$-generation property were
already known (e.g., when $n=1$ and $p=2$ they arise as reflection
groups of the hyperbolic plane with a triangle as a fundamental
domain).  Our construction works for arbitrary $n$ and $p$. For large
$n$ it provides the first examples of non-elementary word hyperbolic
groups possessing the $(n,p)$-generation property.

We construct the groups $G_{n,p}$ as fundamental groups of certain
simplices of groups all of whose local groups are finite $p$-groups.
We use ideas related to simplicial non-positive curvature, developed
by Januszkiewicz and \'Swi\c{a}tkowski in~\cite{JS2}, to show that
these groups are non-elementary word hyperbolic.  The required fixed
point property is obtained using Smith theory and a homological
version of Helly's theorem.

Thus, to complete the proof of Theorem~\ref{thm:fp} and
Corollary~\ref{cor:simp}, it remains to prove Theorems~\ref{thm:comm_quot}
and~\ref{thm:retract_gp}.  This will be done in Sections~\ref{sec:qq}
and~\ref{sec:retract}, respectively.

The validity of  template $\text{\bf NA}_{fp}$ will be established in
Section \ref{sec:triv}; it relies on the Higman Embedding Theorem.  Also
contained in Section~\ref{sec:triv} is Lemma~\ref{lem:simpgp}, which
establishes a triviality property for actions on manifolds.  This is
used both to provide input to the template $\text{\bf NA}_{fp}$ and
to prove Proposition~\ref{prop:simp}.

A version of our paper appeared on the ArXiv preprint server some time
ago.  A number of related papers have since become available,
including \cite{chatterjikassabov,farb,fishersilberman,weinberger}.
(A preliminary version of \cite{farb} predates our work.)
 %Some of the authors were familiar with an early version of
%\cite{farb}, which we gratefully acknowledge.
These papers  consider either actions on manifolds or actions on spaces
equipped with CAT(0)-metrics.  We emphasize that establishing fixed
point properties for actions on arbitrary spaces in $\xac$, or indeed
simplicial actions on simplicial complexes in
$\xac$, is much more difficult.  In particular none
of~\cite{chatterjikassabov,farb,fishersilberman,weinberger} offers an
alternative approach to our main result: the construction of a group
(other than the trivial group) that has the fixed point property
described in Theorem~\ref{thm:fp}.

%%%%%%%%%%%%%%%%%%%%%%%%%%%%%%%%%%%%%%%%%%%%%%%%%%%%%%%
\section{Relatively hyperbolic groups and their common
quotients}\label{sec:qq}
%%%%%%%%%%%%%%%%%%%%%%%%%%%%%%%%

Our purpose in this section is to provide the background we need concerning
relatively hyperbolic groups and their quotients. This will allow us to prove
Propositions~\ref{prop:comm_quot_rh} and ~\ref{prop:comm_quot_hyp} below,
which immediately imply the assertion of Theorem~\ref{thm:comm_quot}.
We adopt the combinatorial approach to
relative hyperbolicity that was developed by Osin in \cite{Osin-RHG}.

Assume that $G$ is a group, $\Hl$ is a fixed collection of proper
subgroups of $G$ (called {\it peripheral subgroups}), and $A$ is a
subset of $G$. The subset  $A$ is called a {\it relative generating
set of $G$} with respect to $\Hl $ if $G$ is generated by $A \cup
\bigcup_{\lambda \in \Lambda} H_\lambda $. In this case $G$ is a
quotient of the free product
%\begin{equation}
$$F=\left( \ast _{\lambda\in \Lambda } H_\lambda  \right) \ast F(A),$$
%\label{F} \end{equation}
where $F(A)$ is the free group with basis $A$. Let $\mathcal R $ be
a subset of $F$ such that the kernel of the natural  epimorphism
$F\twoheadrightarrow G$ is the normal closure of $\mathcal R $ in
the group $F$. In this case we will say that $G$ has the {\it relative
presentation}
\begin{equation}\label{eq:pres_of_G}
\langle A,\; \{H_\lambda\}_{\lambda\in \Lambda}~ \|~ R=1,\,
R\in\mathcal R\rangle .
\end{equation}
If the sets $A$ and $\mathcal R$ are finite, the relative
presentation (\ref{eq:pres_of_G}) is said to be {\it finite}.

Set $\mathcal H=\bigsqcup_{\lambda\in \Lambda} (H_\lambda\setminus
\{ 1\} )$. A finite relative presentation \eqref{eq:pres_of_G} is
said to satisfy a {\it linear relative isoperimetric inequality} if
there exists $C>0$ such that for every word $w$ in the alphabet
$A\cup \mathcal{H}$ (for convenience,  we will further assume that
$A^{-1}=A$) representing the identity in the group $G$, one has
%\begin{equation} \label{prod}
$$w\stackrel{F}{=}\prod\limits_{i=1}^k f_i^{-1}R_i^{\pm 1}f_i,$$
%\end{equation}
with equality in the group $F$, where $R_i\in \mathcal{R}$, $f_i\in
F $, for $i=1, \ldots , k$, and  $k\le C\| w\| $, where $\| w\|$ is
the length of the word $w$.

\begin{definition}\cite{Osin-RHG} The group $G$ is said to be {\it
    relatively hyperbolic} if there is a collection $\Hl$
of proper peripheral subgroups of $G$ such that $G$ admits a finite
relative presentation (\ref{eq:pres_of_G}) satisfying a linear
relative isoperimetric inequality.
\end{definition}

This definition is independent of the choice of the finite
generating set $A$ and the finite set $\mathcal R$ in
(\ref{eq:pres_of_G}) (see \cite{Osin-RHG}).

The definition immediately implies the following basic facts (see
\cite{Osin-RHG}):

\begin{remark} \label{rem:free_prod_rel_hyp}

\rm{(a)} Let $\{H_\lambda\}_{\lambda \in \Lambda}$ be an arbitrary
family of groups. Then the free product $G=\ast_{\lambda\in \Lambda}
H_\lambda$ will be hyperbolic relative to $\Hl$.

{\rm (b)} Any word hyperbolic group (in the sense of Gromov) is
hyperbolic relative to the family $\{\{1\}\}$, where $\{1\}$
denotes the trivial subgroup.
\end{remark}

The following result is our main tool for constructing
common quotients of countable families of relatively
hyperbolic groups. Recall that a group $G$ is said to be {\it
non-elementary} if it does not contain a cyclic subgroup of finite
index.

\begin{thm}\cite[Thm.~1.4]{SQ}\label{thm:comm_quot-SQ}
Any two finitely generated  non-elementary relatively hyperbolic
groups $G_1, G_2$ have a common non-elementary relatively hyperbolic
quotient $H$.
\end{thm}

Consider a sequence of groups $(G_n)_{n\in \mathbb{N}}$ such that
$G_i=G_1/K_i$, $i=2,3,\dots$, for some $K_i \lhd G_1$ and $K_i \le
K_{i+1}$ for all $i \in \N$, $i \ge 2$. The {\it direct limit} of
the sequence $(G_n)_{n\in \mathbb{N}}$ is, by definition, the group
$G_\infty=G_1/K_\infty$ where $K_\infty = \bigcup_{n=2}^\infty K_n$.

\begin{remark} \label{rem:inf_dir_lim} If $G_1$ is finitely generated
  and $G_n$ is infinite for every $n \in \N$,
then $G_\infty$ is also infinite.
\end{remark}

Indeed, suppose that $G_\infty$ is finite, i.e.,
$|G_1:K_\infty|<\infty$. Then $K_\infty$ is finitely generated as a
subgroup of $G_1$, hence there exists $m \in \N$ such that
$K_\infty= K_m$, and $G_\infty=G_m$ is infinite; this is a contradiction.

%\begin{thm}[\cite{SQ}, Thm. 1.4] \label{thm:comm_quot_4_2} Any two finitely generated non-elementary relatively hyperbolic groups $G_1$,$G_2$
%have a common non-elementary relatively hyperbolic quotient Q.
%\end{thm}

\begin{remark} \label{rem:max_norm_sbgp} Any infinite finitely
  generated group $G$ contains a normal subgroup $N$ that is
maximal with respect to the property $|G:N|=\infty$.
\end{remark}

Indeed, let $\mathcal N$ be the set of all normal subgroups of
infinite index in $G$ ordered by inclusion.  Consider a chain
$(M_i)_{i \in I}$ in $\mathcal N$. Set $M=\cup_{i \in I} M_i$; then,
evidently, $M \lhd G$. Now, if $M$ had finite index in $G$, then it
would also be finitely generated. Hence, by the definition of a chain,
there would exist $j \in I$ such that $M=M_j$, which would contradict
the assumption $|G:M_j|=\infty$. Therefore $M \in \mathcal{N}$ is an
upper bound for the chain $(M_i)_{i \in I}$. Consequently, one can
apply Zorn's Lemma to achieve the required maximal element of
$\mathcal N$.

\begin{proposition} \label{prop:comm_quot_rh} Let $\{G_i\}_{i \in \N}$
  be a countable collection of finitely generated
  non-ele\-men\-tary relatively hyperbolic groups and let $C$ be an
  arbitrary countable group. Then there exists a finitely generated
  group $Q$ such that
\begin{itemize}
\item[(i)] $Q$ is a quotient of $G_i$ for every $i \in \N$;
\item[(ii)] $Q$ is a simple group;
\item[(iii)] $Q$ has Kazhdan's property (T);
\item[(iv)] $Q$ contains an isomorphic copy of $C$.
\end{itemize}
\end{proposition}

\begin{proof}%[Proof of Proposition~\label{prop:comm_quot_rh}]
First, embed $C$ into an infinite finitely generated simple group $S$
(see \cite[Ch.~IV, Thm.~3.5]{L-S}). Let $S'$ be a copy of $S$.
Then the group $K=S*S'$ will be non-elementary and hyperbolic
relative to the family consisting of two subgroups $\{S,S'\}$.
Take $G_0$ to be an infinite word hyperbolic group that has property
(T). Then $G_0$ is non-elementary and relatively hyperbolic by
Remark \ref{rem:free_prod_rel_hyp}, hence we can use Theorem
\ref{thm:comm_quot-SQ} to find a non-elementary relatively
hyperbolic group $G(0)$ that is a common quotient of $K$ and $G_0$ (in
particular, $G(0)$ will also be finitely generated).
Now, apply Theorem \ref{thm:comm_quot-SQ} to the groups $G(0)$ and
$G_1$ to obtain their common non-elementary  relatively hyperbolic
quotient $G(1)$. Similarly, define $G(i)$ to be such a quotient for
the groups $G(i-1)$ and $G_i$, $i=2,3,\dots$. Let $G(\infty)$ be the
direct limit of the sequence $\bigl(G(i)\bigr)_{i=0}^\infty$.

The group $G(\infty)$ is finitely generated (as a quotient of $G(0)$)
and infinite (by Remark \ref{rem:inf_dir_lim}), therefore, by Remark
\ref{rem:max_norm_sbgp}, there exists a normal subgroup $N \lhd
G(\infty)$ that is maximal with respect to the property
$|G(\infty):N|=\infty$.  Set $Q=G(\infty)/N$. Then $Q$ is an infinite
group which has no non-trivial normal subgroups of infinite
index. Being a quotient of $G(0)$, makes $Q$ a quotient of $K=S*S'$,
therefore it can not have any proper subgroups of finite index. Thus,
$Q$ is simple. Since the homomorphism $\varphi: S*S' \to Q$ has a
non-trivial image, it must be injective on either $S$ or
$S'$. Therefore $Q$ will contain an isomorphically embedded copy of
$S$, and, consequently, of $C$.

The property (i) for $Q$ follows from the construction.  The
property (iii) holds because $Q$ is a quotient of $G_0$ and since
Kazhdan's property (T) is stable under passing to quotients.
\end{proof}

In the case when one has a collection of word hyperbolic groups (in
the usual, non-relative, sense), one can obtain common quotients with
different properties by using Ol'shanskii's theory of small
cancellation over hyperbolic groups. For example, it is shown in
\cite{Olsh-G-sbgps} that if $g$ is an element of infinite order in a
non-elementary word hyperbolic group $G$, then there exists a number
$n>0$ such that the quotient of $G$ by the normal closure of $g^n$ is
again a non-elementary word hyperbolic group. By harnessing this
result to the procedure for constructing direct limits used in the
proof of Proposition \ref{prop:comm_quot_rh}, we obtain the following
statement, first proved by Osin:

\begin{thm} \cite[Thm.~4.4]{Osin-weakly_amen} \label{thm:per_quot_hyp}
There exists an infinite periodic group $O$, generated by two
elements, such that for every non-elementary word hyperbolic group $H$
there is an epimorphism $\rho: H \twoheadrightarrow O$.
\end{thm}

\begin{proposition} \label{prop:comm_quot_hyp}
%Let $\{G_i\}_{i \in \N}$ be a countable collection of non-ele\-men\-tary word hyperbolic groups.
There exists an infinite finitely generated group $Q$ such that
\begin{itemize}
\item[(a)] $Q$ is a quotient of every non-elementary word hyperbolic group;
\item[(b)] $Q$ is a simple group;
\item[(c)] $Q$ has Kazhdan's property (T);
\item[(d)] $Q$ is periodic.
\end{itemize}
\end{proposition}

\begin{proof} Let $O$ be the group given by Theorem
  \ref{thm:per_quot_hyp}. Since $O$ is finitely generated, it has a
  normal subgroup $N \lhd O$ maximal with respect to the property
  $|O:N|=\infty$ (see Remark \ref{rem:max_norm_sbgp}). Set
  $Q=O/N$. Then $Q$ has no non-trivial normal subgroups of infinite
  index and is a quotient of every non-elementary word hyperbolic
  group; thus $Q$ satisfies (a). In addition, $Q$ is periodic since it
  is a quotient of $O$.

Observe that for an arbitrary integer $k \ge 2$ there exists a
non-elementary word hyperbolic group $H=H(k)$ which does not contain
any normal subgroups of index $k$ (for instance, one can take $H$ to
be the free product of two sufficiently large finite simple groups,
e.g., $H=Alt(k+3)*Alt(k+3)$). Therefore the group $Q$, as a quotient
of $H$, does not contain any normal subgroups of index $k$, for every
$k \ge 2$, hence it is simple. It satisfies Kazhdan's property (T)
because there are non-elementary word hyperbolic groups with (T) and
property (T) is inherited by quotients.
\end{proof}

\begin{remark} The method that we used to obtain simple quotients in
  the proofs of Propositions~\ref{prop:comm_quot_rh} and
\ref{prop:comm_quot_hyp} was highly non-constructive as it relied on
the existence of a maximal normal subgroup of infinite index provided
by Zorn's lemma.  However, one can attain simplicity of the direct
limit in a much more explicit manner, by imposing additional relations
at each step.  For word hyperbolic groups this was done in
\cite[Cor.~2]{Min-p3}.  The latter method for constructing direct
limits of word hyperbolic groups was originally described by
Ol'shanskii in \cite{Olsh-G-sbgps}; it provides significant control
over the resulting limit group. This control allows one to ensure that
the group $Q$ enjoys many properties in addition to the ones listed in
the claim of Proposition \ref{prop:comm_quot_hyp}. For example, in
Proposition~\ref{prop:comm_quot_hyp} one can add that $Q$ has solvable
word and conjugacy problems.
\end{remark}

%%%%%%%%%%%%%%%%%%%%%%%%%%%%%%%%%%%%%%%%%%%%%%%%%%%%%%%%%%%%%%%%
\section{Simplices of finite $p$-groups with non-elementary word
  hyperbolic direct limits}\label{sec:retract}
%%%%%%%%%%%%%%%%%%%%%%%%%%%%%%%%%%%%%%%%%%%%%%%%%%%%%%%%%%%%%%%%

Theorem \ref{thm:retract_gp} is an immediate consequence of the
following  two results, whose proof is the object of this section.

\begin{thm}\label{thm:monstrous}
For every prime number $p$ and integer $n\ge1$
there is a non-elementary word hyperbolic group $G$
generated by a set $S$ of cardinality $n+2$
such that the subgroup of $G$ generated by each proper
subset of $S$ is a finite $p$-group.
\end{thm}

\begin{thm}\label{thm:helly}
Let $p$ be a prime number.
Suppose that a group $G$ has a generating set $S$ of cardinality $n+2$,
such that the subgroup generated by each proper subset of $S$ is a
finite $p$-group.  Then
for any action of $G$ on a Hausdorff mod-$p$ acyclic space $X$
of covering dimension less than or equal to $n$, the global fixed
point set is mod-$p$ acyclic.
\end{thm}

We prove Theorem \ref{thm:monstrous} by constructing each of the
desired groups as the fundamental group (equivalently, the direct
limit) of a certain $(n+1)$-dimensional simplex of finite $p$-groups.
The construction of the local groups in each simplex of groups, for
fixed $p$, proceeds by induction on $n$.  Provided that $m\le n$, the
groups that are assigned to each codimension $m$ face of the
$(n+1)$-simplex will depend only on $m$, up to isomorphism.  The
codimension zero face, i.e., the whole $(n+1)$-simplex itself, will be
assigned the trivial group 1, and each codimension one simplex will be
assigned a cyclic group of order $p$. As part of the inductive step,
we will show that the fundamental group of the constructed
$(n+1)$-simplex of groups maps onto a $p$-group in such a way that
each local group maps injectively.  This quotient $p$-group will be
the group used as each vertex group in the $(n+2)$-simplex of groups.

The idea that drives our construction consists of
requiring and exploiting existence of certain retraction homomorphisms
between the local groups of the complexes of groups involved.
We develop this approach in Subsections 3.2--3.4 below,
after recalling in Subsection 3.1 some basic notions and facts
related to complexes of groups.

Each simplex of groups that we construct will be developable.
Associated to any developable $n$-simplex of groups $\bfg$, there is a
simplicial complex $X$ on which the fundamental group $G$ of $\bfg$
acts with an $n$-simplex as strict fundamental domain.  If the local
groups of $\bfg$ are all finite, the corresponding action is
proper. Thus we may show that the group $G$ is word hyperbolic by
showing that the associated simplicial complex $X$ is Gromov
hyperbolic.  We show that $X$ is indeed hyperbolic by verifying that
it satisfies a combinatorial criterion for the hyperbolicity of a
simplicial complex related to the idea of simplicial non-positive curvature
developed in \cite{JS2}.  More precisely, we show that $X$ is
8-systolic, and hence hyperbolic.  This is the content of
Subsection~3.5.

From the perspective of the subject of simplicial non-positive curvature,
Subsections 3.1-3.5 may be viewed as providing an alternative to the
construction from~\cite{JS2} of numerous examples of $k$-systolic groups
and spaces, for arbitrary $k$ and in arbitrary dimension.
The resulting groups are different from those obtained in \cite{JS2}.

From the algebraic perspective, this construction provides new operations
of product type for groups, the so called $n$-{\it retra-products},
which interpolate between the direct product and the free product.
These operations can be further generalized in the spirit of graph products.
We think the groups obtained this way deserve to be studied.
Such groups fall in the class
of systolic, or even 8-systolic groups, and thus share various
exotic properties of the latter, as established in
\cite{JS2, JS3, Os1, Os2}.
In a future work we plan to show that
$n$-retra-products of finite groups, for sufficiently large $n$,
are residually finite.

The last subsection of this section, Subsection 3.6, contains the
proof of Theorem \ref{thm:helly}.  This proof uses a result from Smith
theory concerning mod-$p$ cohomology of the fixed point set of a
finite $p$-group action.  It also uses a homological version of
Helly's theorem for mod-$p$ acyclic subsets.

\begin{remark}
It is because we need to apply Smith theory that the groups previously
constructed in \cite{JS2} are unsuitable for our purposes.  The groups
constructed in \cite{JS2} include fundamental groups of simplices of
\emph{finite} groups which are non-elementary word hyperbolic, but the
finite groups occurring in \cite{JS2} are not of prime power order.
\end{remark}

\begin{remark} A simplified form of the arguments in Subsection 3.6
shows that the fundamental group
of any $(n+1)$-simplex of finite groups cannot act without a global fixed
point by isometries on any complete CAT(0) space of covering dimension
at most $n$.  Indeed, the Helly-type argument we use
goes through almost unchanged,
while the fact that the fixed point set for a finite group of
isometries of a complete CAT(0) space is contractible replaces the
appeal to Smith theory.  Versions of this argument  have appeared
previously  in
work of Barnhill~\cite{barnhill}, Bridson (unpublished), and
Farb~\cite{farb}.
\end{remark}

\subsection{Strict complexes of groups}

We recall some basic notions and facts related to strict complexes
of groups. The main reference is Bridson and Haefliger~\cite{BH},
where these objects are called simple complexes of groups.

A simplicial complex $K$ gives rise to two categories: the category
${\mathcal Q}_K$ of non-empty simplices of $K$ with inclusions as
morphisms, and the extended category ${\mathcal Q}^+_K$ of simplices
of $K$ including the empty set $\emptyset$ as the unique
$(-1)$--simplex.  In addition to the morphisms from ${\mathcal Q}_K$,
the category ${\mathcal Q}^+_K$ has one morphism from $\emptyset$ to
$\sigma$ for each nonempty simplex $\sigma$ of $K$.  A {\it strict
complex of groups} $\bfg$ consists of a simplicial complex $|\bfg|$
(called {\it the underlying complex of $\bfg$}), together with a
contravariant functor $\bfg$ from ${\mathcal Q}_{|\bfg|}$ to the
category of groups and embeddings.  A strict complex of groups is {\it
developable} if the functor $\bfg$ extends to a contravariant functor
$\bfg^+$ from the category ${\mathcal Q}^+_{|\bfg|}$ to the category
of groups and embeddings.  Given an extension $\bfg^+$ of $\bfg$, we
will denote by $G$ the group $\bfg^+(\emptyset)$. For simplices
$\tau\subset\sigma$ (allowing $\tau=\emptyset$), we will view the
group $\bfg(\sigma)$ as subgroup in the group $\bfg(\tau)$. We will be
interested only in extensions that are {\it surjective}, i.e.  such
that the group $G=\bfg^+(\emptyset)$ is generated by the union of its
subgroups $\bfg(\sigma)$ with $\sigma\ne\emptyset$.

We call any surjective extension $\bfg^+$ of $\bfg$ an {\it extended
complex of groups}.  We view the collection of all possible surjective
extensions of $\bfg$ to ${\mathcal Q}^+_{|\bfg|}$ also as a category,
which we denote by ${\rm Ext}_\bfg$. We take as morphisms of ${\rm
Ext}_\bfg$ the natural transformations from $\bfg^+$ to ${\bfg^+}'$
which extend the identical natural transformation of $\bfg$.  (Note
that, given extensions $\bfg^+$ and ${\bfg^+}'$, there may be no
morphism between them, and if there is one then, by surjectivity, it
is unique; moreover, the homomorphism from $G=\bfg^+(\emptyset)$ to
$G'={\bfg^+}'(\emptyset)$ induced by a morphism is {\sl not} required
to be an embedding, although it is required to be a group
homomorphism.)  If $\bfg$ is developable, the category ${\rm
Ext}_\bfg$ has an initial object $\bfg^+_{dir}$, in which the group
$\bfg^+_{dir}(\emptyset)$ is just the direct limit of the functor
$\bfg$ (for brevity, we often denote this direct limit group $\tilG$).
Thus for any extension $\bfg^+$ of $\bfg$, there is a unique group
homomorphism from $\tilG$ to $G$ extending the identity map on each
$\bfg(\sigma)$ for $\sigma\in {\mathcal Q}_{|\bfg|}$.  In the cases
that will be considered below, the simplicial complex $|\bfg|$ is
simply connected, which implies that $\tilG$ coincides with what is
known as the fundamental group of the complex of groups $\bfg$.  (In
fact, $|\bfg|$ is contractible in the cases considered below.)

For an extended complex of groups $\bfg^+$ we consider a space
$d\bfg^+$ with an action of $G=\bfg^+(\emptyset)$, the {\it
development} of $\bfg^+$, given by
$$
d\bfg^+= |{\bfg}|\times G/_\sim,
$$
where the equivalence relation $\sim$ is given by $(p,g)\sim(q,h)$
iff $p=q$ and there exists $\sigma\in {\mathcal Q}_{|\bfg|}$
so that $p\in \sigma$ and
$g^{-1}h\in \bfg(\sigma)$.  It suffices to take $\sigma$ to be the
minimal simplex containing $p$.  The $G$-action is given by $g[p,h]=
[p,gh]$.  The quotient by the action of $G$ is
(canonically isomorphic with) $|\bfg|$, and the
subcomplex
$$
[|\bfg|,1]=\{ [(p,1)]:p\in|\bfg| \},
$$
(where $[(p,1)]$ is the equivalence class of $(p,1)$ under $\sim$)
is a strict fundamental domain for the action (in
the sense that every $G$-orbit intersects $[|\bfg|,1]$ in exactly one
point).  The space $d\bfg^+$ is a multi-simplicial complex, and the
(pointwise and setwise) stabilizer of the simplex $[\sigma,g]$ is
the subgroup $g\bfg(\sigma)g^{-1}$. In the cases considered below,
developments will be true simplicial complexes.

A {\it morphism} $\varphi$ from a strict complex of groups $\bfg$
to a group $H$ is a compatible collection of homomorphisms
$\varphi_\sigma:\bfg(\sigma)\to H, \sigma\in{\mathcal Q}_{|\bfg|}$
(in general not necessarily injective).
{\it Compatibility} means that we have equalities
$\varphi_\sigma=\varphi_\tau\circ i_{\sigma\tau}$
for any $\tau\subset\sigma$, where $i_{\sigma\tau}$
is the inclusion of $\bfg(\sigma)$ in $\bfg(\tau)$.
For example, a collection of inclusions $\bfg(\sigma)\to\bfg^+(\emptyset)$
is a morphism $\bfg\to\bfg^+(\emptyset)$.
A morphism $\varphi:\bfg\to H$ is {\it locally injective}
if all the homomorphisms $\varphi_\sigma$ are injective.

Suppose we are given an action of a group $H$ on a simplicial
complex $X$, by simplicial automorphisms, and suppose this action is
{\it without inversions}, i.e., if $g\in H$ fixes a simplex of $X$,
it also fixes all vertices of this simplex. Suppose also that the
action has a strict fundamental domain $D$ which is a subcomplex of
$X$. Clearly, $D$ is then isomorphic to the quotient complex
$H\backslash X$. Such an action determines the {\it extended
associated complex of groups} $\bfg^+$, with the underlying complex
$|\bfg|=D$, with local groups
$\bfg(\sigma):=Stab(\sigma,H)$ for $\sigma\subset D$, and with
$\bfg^+(\emptyset):=H$. The morphisms in $\bfg$ are the natural inclusions. It
turns out that in this situation the development $d\bfg^+$ is $H$-equivariantly
isomorphic with $X$.

%**********************************
\subsection{Higher retractibility}
%**********************************

Now we pass to a less standard part of the exposition. We begin by
describing a class of simplicial complexes, called {\it blocks},
that will serve as  the underlying complexes of the complexes of groups
involved in our construction. Then we discuss various requirements
on the corresponding complexes of groups. Some part of this material
is parallel to that in Sections 4 and 5 of~\cite{JS1}, where
retractibility and extra retractibility stand for what we call in
this paper 1-retractibility and 2-retractibility, respectively.

\begin{definition}[Block]\label{1} A simplicial complex $K$ of dimension
$n$ is a  \emph{chamber complex} if each of its simplices is a face
of an $n$-simplex of $K$. Top dimensional simplices are then called
\emph{chambers} of $K$. A chamber complex $K$ is \emph{gallery
connected} if each pair of its chambers is connected by a sequence of
chambers in which any two consecutive chambers share a face of
codimension 1. A chamber complex is \emph{normal} if it is gallery
connected and all of its links (which are also chamber complexes)
are gallery connected. The  \emph{boundary} of a chamber complex
$K$, denoted $\partial K$, is the subcomplex of $K$ consisting of
all those faces of codimension 1 that are contained in precisely
one chamber. A \emph{block} is a normal chamber complex with
nonempty boundary. The \emph{sides} of a block $B$ are the faces of
codimension 1 contained in $\partial B$. We denote the set of all
sides of $B$ by ${\mathcal S}_B$.
\end{definition}

Note that links $B_\sigma$ of a block at faces
$\sigma\subset\partial B$ are also blocks, and that
$\partial(B_\sigma)=(\partial B)_\sigma$.

\begin{definition}[Normal block of groups]\label{2} A \emph{normal block of
groups} over a block $B$ is a strict complex of groups $\hbox{\bf
G}$ with $|\bfg|=B$ satisfying the following two conditions:
\begin{itemize}
  \item[(1)] $\hbox{\bf G}$ is \emph{boundary supported},
  i.e. $\hbox{\bf G}(\sigma)=1$
for each $\sigma$ not contained in $\partial B$;
  \item[(2)] $\hbox{\bf G}$ is \emph{locally $\mathcal S$-surjective},
  i.e., every group $\hbox{\bf G}(\sigma)$
is generated by the union $\bigcup \{ \hbox{\bf G}(s):s\in{\mathcal
S}_B, \sigma\subset s \}$, where we use the convention that the
empty set generates the trivial group 1.
\end{itemize}

\noindent
An \emph{extended normal block of groups} is an extension $\hbox{\bf
G}^+$ of a normal block of groups $\bfg$
such that the associated morphism
$\varphi:\hbox{\bf G}\to \hbox{\bf G}(\emptyset)$  is
$\mathcal S$-\emph{surjective},
i.e. $\hbox{\bf G}(\emptyset)$ is generated by the
union $\bigcup \{ \hbox{\bf G}(s):s\in{\mathcal S}_B \}$.
(To simplify notation, we write $\bfg(\sigma)$ instead of $\bfg^+(\sigma)$
to denote the corresponding groups of $\bfg^+$.)
\end{definition}

Given an extended normal block of groups $\hbox{\bf G}^+$ over $B$,
its development $ d\bfg^+$ is tesselated by copies of $B$. More
precisely, $ d\bfg^+$ is the union of the subcomplexes of the form
$[B,g]$, with $g\in\hbox{\bf G}(\emptyset)$, which do not intersect
each other except at their boundaries, and which we view as tiles of
the tesselation. Moreover, the action of $\hbox{\bf G}(\emptyset)$
on $ d\bfg^+$ is simply transitive on these tiles. By $\mathcal S$-surjectivity
of $\varphi$, $ d\bfg^+$ is a normal chamber complex. If $|\hbox{\bf
G}(s)|>1$ for all $s\in{\mathcal S}_B$ then the chamber complex $
d\bfg^+$ has empty boundary. If $B$  is a pseudo-manifold and
$|\hbox{\bf G}(s)|\le2$ for all $s\in{\mathcal S}_B$, then $
d\bfg^+$ is a pseudo-manifold.

\begin{definition}[1-retractibility]\label{3} An extended normal block of
groups $\hbox{\bf G}^+$ is 1-\emph{retrac\-tible} if for every
$\sigma\subset |\bfg|$ there is a homomorphism $r_\sigma:\hbox{\bf
G}(\emptyset)\to \hbox{\bf G}(\sigma)$ such that
$r_\sigma|_{{\bf G}(s)}=id_{{\bf G}(s)}$ for
$s\in{\mathcal S}_{|\bfg|}, s\supset\sigma$, and
$r_\sigma|_{{\bf G}(s)}=1$ otherwise.
\end{definition}

In the next two lemmas we present properties that immediately follow
from 1-retracti\-bi\-lity. We omit the straightforward proofs.

\begin{lemma}\strut
\begin{itemize}
  \item[(1)] The homomorphisms $r_\sigma$, if they exist, are unique.
  \item[(2)] Let $\varphi_\sigma:\hbox{\bf G}(\sigma)\to
  \hbox{\bf G}(\emptyset)$ be the
homomorphisms of the morphism $\varphi:\hbox{\bf G}\to \hbox{\bf
G}(\emptyset)$. Then for each $\sigma$ we have
$r_\sigma\varphi_\sigma=id_{{\bf G}(\sigma)}$. Thus $r_\sigma$
is a retraction onto the subgroup $\hbox{\bf G}(\sigma)<\hbox{\bf
G}(\emptyset)$.
  \item[(3)] The inclusion homomorphisms $\varphi_{\tau\sigma}:\hbox{\bf G}(\tau)\to
\hbox{\bf G}(\sigma)$, for $\sigma\subset\tau$, occurring as the
structure homomorphisms of $\hbox{\bf G}$, satisfy
$\varphi_{\tau\sigma}=r_\sigma\varphi_\tau$.
\end{itemize}
\end{lemma}

Motivated by property (3) above, we define homomorphisms
$r_{\rho\tau}:\hbox{\bf G}(\rho)\to \hbox{\bf G}(\tau)$, for any
simplices $\rho,\tau$ of $|\bfg|$, including $\emptyset$, by putting
$r_{\rho\tau}:=r_\tau\varphi_\rho$.

\begin{lemma}\label{5}  Each of the homomorphisms $r_{\rho\tau}$ is
uniquely determined by the following two requirements:
\begin{itemize}
  \item[(1)] $r_{\rho\tau}|_{{\bf G}(s)}=id_{{\bf G}(s)}$ for $s\in{\mathcal S}_{|\bfg|},
s\supset\rho, s\supset\tau$;
  \item[(2)] $r_{\rho\tau}|_{{\bf G}(s)}=1$ otherwise (i.e. for
$s\in{\mathcal S}_{|\bfg|}, s\supset\rho, s\hbox{ not containing
}\tau$).
\end{itemize}
 In particular, we have
$r_{\sigma\emptyset}=\varphi_\sigma$,
$r_{\emptyset\sigma}=r_\sigma$, and
$r_{\tau\sigma}=\varphi_{\tau\sigma}$ whenever $\tau\supset\sigma$.
Moreover, if $\tau\supset\sigma$ then $r_{\tau\sigma}$ is a
retraction (left inverse) for the inclusion $\varphi_{\tau\sigma}$.
\end{lemma}

To define higher retractibility properties for an extended normal
block of groups $\hbox{\bf G}^+$ we need first to introduce certain
new blocks of groups called \emph{unfoldings} of $\hbox{\bf G}^+$ at
the boundary simplices $\sigma\subset\partial |\bfg|$.

\begin{definition}[Unfolding of $\hbox{\bf G}^+$ at $\sigma$]\label{6} Let
$\hbox{\bf G}^+$ be a 1-retractible extended normal block of groups.
Let $\sigma\subset\partial |\bfg|$ be a simplex, and denote by
$d_\sigma \hbox{\bf G}(\emptyset)$ the kernel of the retraction
homomorphism $r_\sigma:\hbox{\bf G}(\emptyset)\to \hbox{\bf
G}(\sigma)$. The \emph{unfolding} of $\hbox{\bf G}^+$ at $\sigma$,
denoted $ d_\sigma\bfg $, is the complex of groups associated to the
action of the group $ d_\sigma\bfg (\emptyset)$ on the development $
d\bfg^+$. The \emph{extended unfolding} $d_\sigma\bfg^+ $ is the
same complex of groups equipped with the canonical morphism to the
group $ d_\sigma\bfg (\emptyset)$.
\end{definition}

Define a subcomplex $d_\sigma|\bfg|\subset  d\bfg^+$ by
$d_\sigma|\bfg|:=\bigcup\{ [B,g] : g\in\hbox{\bf G}(\sigma) \}$. We
will show that the above defined unfolding $d_\sigma\bfg^+ $ is an
extended normal block of groups over $d_\sigma|\bfg|$. This will be
done in a series of lemmas, in which we describe the structure of
$d_\sigma|\bfg|$ and $d_\sigma\bfg^+ $ in detail.

\begin{lemma}\label{7} $d_\sigma|\bfg|$ is a strict fundamental domain
for the action of the group $ d_\sigma\bfg (\emptyset)$ on the
development $ d\bfg^+$. In particular,
$|d_\sigma\bfg|=d_\sigma|\bfg|$.
\end{lemma}

\begin{proof} We need to show that the restriction to
$d_\sigma|\bfg|$ of the quotient map $q_\sigma: d\bfg^+\to
d_\sigma\bfg (\emptyset)\backslash d\bfg^+$ is a bijection. This
follows by observing that the map $j_\sigma:  d_\sigma\bfg
(\emptyset)\backslash d\bfg^+\to d_\sigma|\bfg|$ defined by
$j_\sigma(  d_\sigma\bfg (\emptyset)\cdot[p,g])=[p,r_\sigma(g)]$ is
the inverse of $q_\sigma|_{d_\sigma|\bfg|}$. \end{proof}

To proceed with describing $d_\sigma|\bfg|$, we need to define links
for blocks of groups. This notion will  also be useful in our later
considerations.

\begin{definition}[Link of a block of groups]\label{8} Let $\bfg$ be a
normal block of groups and let $\sigma$  be a simplex of $|\bfg|$.
The \emph{link} of $ \bfg^+$ at $\sigma$ is an extended normal block
of groups $\hbox{\bf G}^+_\sigma$ over the link $|\bfg|_\sigma$
given by $\hbox{\bf G}^+_\sigma(\tau):=\hbox{\bf G}(\tau*\sigma)$ for
all $\tau\subset |\bfg|_\sigma$, including the empty set $\emptyset$
(with the convention that $\emptyset*\sigma=\sigma$).
\end{definition}

We skip the straightforward argument for showing that the above
defined extended complex of groups is an extended normal block of
groups.

The next lemma describes the links of the complex $d_\sigma|\bfg|$.
We omit its straightforward proof.
In this lemma, and in the remaining part of this section,
we will denote by
$\sigma-\tau$  the face of $\sigma$
spanned by the vertices of $\sigma$ not contained in $\tau$.

\begin{lemma}\label{9} Let $[\tau,g]$ be a simplex of $d_\sigma|\bfg|$,
where $\tau\subset |\bfg|$ and $g\in\hbox{\bf G}(\sigma)$. For any
simplex $\rho\subset |\bfg|_\tau$ let $d_\rho |\bfg_\tau|$ be the
strict fundamental domain for the action of the group
$d_\rho\hbox{\bf G}_\tau(\emptyset)$ on the development $d\hbox{\bf
G}_\tau^+$. Then the link of $d_\sigma|\bfg|$ at $[\tau,g]$ has one
of the following two forms depending on $\tau$:
\begin{itemize}
  \item[(1)] $(d_\sigma|\bfg|)_{[\tau,g]}\cong d_{\sigma-\tau}|\bfg_\tau|$
if $\sigma$  and $\tau$ span
a simplex of $|\bfg|$,
where we use convention that $d_\emptyset|\bfg_\tau|=d\hbox{\bf
G}_\tau^+$;
  \item[(2)] $(d_\sigma|\bfg|)_{[\tau,g]}\cong |\bfg|_\tau$ otherwise.
\end{itemize}\end{lemma}

Lemma~\ref{9} easily implies the following corollary. The proof of
part (1) uses induction on the dimension of $B$ and $\mathcal
S$-surjectivity of the extending morphism; we omit the details.

%\bigskip
%\eject

\begin{coroll}\label{10}\strut
\begin{itemize}\item[(1)] $d_\sigma|\bfg|$ is a normal chamber complex.
\item[(2)]   The boundary $\partial(d_\sigma|\bfg|)$
is the subcomplex of $d_\sigma|\bfg|$
consisting of the simplices of form $[\rho,g]$ for all
$\rho\subset\partial |\bfg| $ not containing $\sigma$ and for all
$g\in\bfg(\sigma)$. In particular, the set of sides of
$d_\sigma|\bfg|$ is the set
$${\mathcal S}_{d_\sigma|\bfg|}=\{ [s,g]
: s\in{\mathcal S}_{|\bfg|}, s \hbox{ does not contain } \sigma,
g\in\bfg(\sigma)  \}
.$$\end{itemize}\end{coroll}

For a subgroup $H<G$ and an element $g\in G$, we denote by $H^g$ the
conjugation $gHg^{-1}$.  The next lemma describes the local groups
of the unfolding $ d_\sigma\bfg $.

\begin{lemma}\label{11}  Let $ d_\sigma\bfg $ be the unfolding of $\bfg$ and
let $[\tau,g]$ be a simplex of $d_\sigma|\bfg|$ (with $\tau\subset
|\bfg|$ and $g\in\bfg_\sigma$). Then
$$
 d_\sigma\bfg ([\tau,g])=[\ker(r_{\tau\sigma}:
 \bfg(\tau)\to\bfg(\sigma))]^g< d_\sigma\bfg (\emptyset),
$$
and consequently
$$
 d_\sigma\bfg ([\tau,g])=g\cdot \langle \bigcup\bfg(s):s\in{\mathcal S}_{|\bfg|},
 s\supset\tau,
s \hbox{ does not contain }\sigma \rangle \cdot g^{-1}.
$$
\end{lemma}

In view of Corollary~\ref{10}(2), Lemma~\ref{11} implies the
following.

\begin{coroll}\label{12} For simplices $[\tau,g]$ of $d_\sigma|\bfg|$
not contained in the boundary $\partial(d_\sigma|\bfg|)$ we have $
d_\sigma\bfg ([\tau,g])=1$.\end{coroll}

Another consequence of Lemma~\ref{11}, which will be useful later,
is the following.

\begin{coroll}\label{13}  Let $\bfg^+$ be an extended normal block of
groups, and let $\sigma$ be a simplex of $|\bfg|$. Then for any
simplex $[\tau,g]\subset d_\sigma|\bfg|$, with $\tau\subset|\bfg|$
and $g\in\bfg(\sigma)$, we have $[d_\sigma\bfg]_{[\tau,g]} \cong
d_{\sigma-\tau}\bfg_\tau^+$, where $\cong$ denotes an isomorphism of
extended complexes of groups, and where $d_\emptyset\bfg_\tau^+$
denotes here the trivial strict complex of groups over $d\bfg_\tau^+$
(i.e. all of the local groups are trivial).
\end{coroll}

As a consequence of the results above, from Lemma~\ref{7} to
Corollary~\ref{12}, we obtain the following.

\begin{coroll}\label{14}  Given an extended normal block of groups
$\bfg^+$, each of its unfoldings $d_\sigma\bfg^+ $ is an extended
normal block of groups.
\end{coroll}

We are now in a position to define recursively higher
retractibilities.

\begin{definition}[$n$-retractibility]\label{15} Let $n$ be a natural number.
An extended normal block of groups $\bfg^+$ is \emph{
$(n+1)$-retractible} if it is 1-retractible, and for every simplex
$\sigma\subset\partial |\bfg|$ the unfolding $d_\sigma\bfg^+ $ is
$n$-retractible.
\end{definition}

\begin{example}[$n$-retractible 1-simplex of groups]\label{16} Consider
the extended complex of groups $\bfg^+$ with $|\bfg|$ equal to a
1-simplex, with the vertex groups $\bfg(v)$ cyclic of order two, and
with $\bfg(\emptyset)$ dihedral of order $2k$, where generators of
the vertex groups correspond to standard generators of the dihedral
group.  This complex of groups is clearly an extended normal block
of groups. Moreover, it is $n$-retractible but not
$(n+1)$-retractible in the case when $k=2^n(2m+1)$ for some $m$.
\end{example}

\begin{remark}\label{17}\strut
\begin{itemize}\item[(1)] Note that an $n$-retrac\-tible extended normal block of groups is
always $k$-retrac\-tible for each $k\leq n$.
\item[(2)] All links in a 1-retractible normal block of groups are 1-retrac\-tible.
\item[(3)] From Corollary \ref{13} and the above remark (2)
one can easily deduce using induction on $n$ that
if $\bfg^+$ is $n$-retractible and
$\sigma$ is a simplex of $|\bfg|$, then $\bfg_\sigma^+$ is also
$n$-retractible.\end{itemize}
\end{remark}

%***********************************************************
\subsection{$n$-retractible extensions}\label{1.3}
%**********************************************************

Our next objective is to establish results that give partial converses to
property in Remark~\ref{17}(3). These will allow us to pass up one
dimension in our recursive construction of $n$-retractible
simplices of groups, in the next subsection.

\begin{proposition}\label{18} Let $\bfg$ be a normal block of groups.
If for some natural number
$n$ all the links in $\bfg$ (as extended complexes of groups)
are $n$-retractible then:
\begin{itemize}
\item[(1)] $\bfg$ is developable, and
\item[(2)] the extension $\bfg^+_{dir}$ of $\bfg$ (in which
${\bfg}(\emptyset)$ coincides with the direct limit $\tilG$ of $\bfg$)
is $n$-retractible.
\end{itemize}
\end{proposition}

\begin{proof}
To prove part (1), we need to show that $\bfg$ admits
a locally injective morphism $\psi:\bfg\to H$ to some group $H$.
To find $\psi$, for each $\sigma\subset|\bfg|$ we construct
a morphism $\bar r_\sigma:\bfg\to\bfg(\sigma)$
which is identical on the group $\bfg(\sigma)$ of $\bfg$.
We then take as $H$ the direct product
$H=\oplus\{ \bfg(\sigma):\sigma\subset|\bfg| \}$
and as $\psi$ the diagonal morphism
$\psi=\oplus\{ \bar r_\sigma:\sigma\subset|\bfg| \}$,
which is then clearly locally injective.

Fix a simplex $\sigma\in |\bfg|$.
To get a morphism $\bar r_\sigma$ as above,
we will construct an appropriate compatible collection
of homomorphisms
$\bar r_{\eta\sigma}:\bfg(\eta)\to\bfg(\sigma)$,
for all $\eta\subset|\bfg|$, such that
$\bar r_{\sigma\sigma}=id_{\bfg(\sigma)}$.
To do this,
for any simplex $\rho\subset|\bfg|$ consider the set of sides
${\mathcal S}_\rho=\{ s\in{\mathcal S}_{|\bfg|} :\rho\subset s \}$.
Fix a  simplex $\eta\subset |\bfg|$.
If ${\mathcal S}_\eta\cap{\mathcal S}_\sigma=\emptyset$,
put $\bar r_{\eta\sigma}$ to be the trivial homomorphism.
Otherwise, consider the simplex
$\tau=\cap\{ s:s\in{\mathcal S}_\eta\cap{\mathcal S}_\sigma$ \}.
Clearly, we have then $\eta\subset\tau$ and $\sigma\subset\tau$.
In particular, we have the inclusion homomorphism
$i_{\tau\sigma}:\bfg(\tau)\to\bfg(\sigma)$,
which is identical on the subgroups
$\bfg(s):s\in{\mathcal S}_\eta\cap{\mathcal S}_\sigma$.

Recall that we denote by $\tau-\eta$ the face of $\tau$
spanned by all vertices not contained in $\eta$.
Since $\eta\subset\tau$, the
groups $\bfg(\eta)$ and $\bfg(\tau)$ coincide with the link groups
$\bfg^+_\eta(\emptyset)$ and $\bfg^+_\eta(\tau-\eta)$, respectively.
Since the link $\bfg_\eta^+$ is 1-retractible,
we have the retraction
$$
r_{\tau-\eta}:
\bfg^+_\eta(\emptyset)\to\bfg^+_\eta(\tau-\eta)
$$
such that
$r_{\tau-\eta}|_{\bfg^+_\eta(s)}=id_{\bfg^+_\eta(s)}$
for $s\in{\mathcal S}_{|\bfg|_\eta}, s\supset\tau-\eta$
and $r_{\tau-\eta}|_{\bfg^+_\eta(s)}=1$ otherwise
(cf. Definition \ref{3}).

Put
$\bar r_{\eta\sigma}:=i_{\tau\sigma}\circ r_{\tau-\eta}$.
We claim that $\bar r_{\eta\sigma}$ satisfies the assertions
(1) and (2) of Lemma \ref{5}, when substituted for $r_{\eta\sigma}$.
This follows from
the identification of the groups $\bfg(s), s\in{\mathcal S}_{|\bfg|}$
with the groups $\bfg^+_\eta(s-\eta)$, for $\eta\subset s$,
and from the fact that $s-\eta\in{\mathcal S}_{|\bfg|_\eta}$
(because $\partial(|\bfg|_\eta)=(\partial|\bfg|)_\eta$).

Now, we need to check the compatibility condition
$\bar r_{\eta_2\sigma}=\bar r_{\eta_1\sigma}\circ i_{\eta_2\eta_1}$
for all simplices $\eta_1\subset\eta_2$ in $|\bfg|$.
This follows from the coincidence of the maps on both sides
of the equality on the generating set
$\cup\{ \bfg(s):s\in{\mathcal S}_{|\bfg|},  s\supset\eta_2,  \}$
of $\bfg(\eta_2)$.
This coincidence is a fairly direct consequence of
assertions (1) and (2) of Lemma \ref{5}, satisfied
by the maps $\bar r_{\eta_1\sigma}$ and $\bar r_{\eta_2\sigma}$.

Finally, assertion (1) of Lemma \ref{5} clearly implies
that $\bar r_{\sigma\sigma}=id_{\bfg(\sigma)}$,
which concludes the proof of developability of $\bfg$.

We now turn to proving part (2).
To deal with the case $n=1$ we need to construct
the map $r_\sigma:\tilG\rightarrow \bfg(\sigma)$ as required in
Definition~\ref{3}, for any $\sigma\subset|\bfg|$.

Consider the maps $\bar r_{\eta\sigma}:\eta\subset|\bfg|$
constructed in the proof of part (1),
and the morphism $\bar r_\sigma:\bfg\to\bfg(\sigma)$
given by these maps. Let $r_\sigma:\tilG\to\bfg(\sigma)$
be the homomorphism induced by this morphism.
The requirements of Definition \ref{3} for $r_\sigma$ follow then easily
from the assertions of Lemma \ref{5} satisfied by the
maps $\bar r_{\eta\sigma}$ (we skip the straightforward details).
Thus 1-retractibility of links of $\bfg$ implies 1-retractibility
of $\bfg^+_{dir}$.

Now suppose that $n>1$. If links in $\bfg$ are $n$-retractible, it
follows that links in $d_\sigma\bfg$ are $(n-1)$-retractible for all
$\sigma$.  By induction, it follows that the unfoldings
$d_\sigma\bfg^+_{dir}$ are $(n-1)$-retractible, and so
$\bfg^+_{dir}$ is $n$-retractible. \end{proof}

\begin{example}[$n$-retractible 2-simplex of groups]\label{19} Consider the
triangle Coxeter group $W_{n,m}$ given by
$$
W_{n,m}=\langle s_1,s_2,s_3|s_i^2, (s_is_j)^k \hbox{ for } i\ne j
\rangle,
$$
where $k=2^n(2m+1)$. The Coxeter complex for this group is a
triangulation of the (hyperbolic) plane on which the group acts by
simplicial automorphisms, simply transitively on 2-simplices. Let
$\bfg^+$ be the extended 2-simplex of groups associated to this
action. Links in $\bfg$ at vertices are then isomorphic to the 1-simplex
of groups from Example 3.18. Thus, in view of Proposition~\ref{18} and
Remark~\ref{17}(3), $\bfg^+$ is $n$-retractible, but not
$(n+1)$-retractible.
\end{example}

A more subtle way of getting $n$-retractible extensions is given in
the following theorem, which will be used directly, as a recursive
step, in our construction in Subsection~\ref{1.4}.

\begin{thm}\label{20} Let $\bfg$ be a normal block of groups in
which all links are $n$-retractible.  Then there exists an extension
$\bfg^+_\min$ of $\bfg$ that has the following properties.

\begin{itemize}
\item[(1)] $\bfg^+_\min$ is the minimal $n$-retractible extension
of $\bfg$ in the following sense: if $\bfg^+$ is any $n$-retractible
extension of $\bfg$ then there is a unique morphism of extended
complexes of groups $\bfg^+\to \bfg^+_\min$ which extends the
identity on $\bfg$ (i.e. a morphism in the category ${\rm Ext}_\bfg$). \pra
\item[(2)] If all $\bfg(\sigma)$ are finite, so is $\bfg_\min(\emptyset)$.
\pra
\item[(3)] If all $\bfg(\sigma)$ are $p$-groups of bounded exponent,
so is $\bfg_\min(\emptyset)$. \pra
\item[(4)] If all $\bfg(\sigma)$ are soluble groups,
%with soluble length $\le d$ for some fixed $d$,
then so is $\bfg_\min(\emptyset)$.
\end{itemize}
\end{thm}

\begin{remark} Property 1 means that the extended complex of
groups $\bfg^+_\min$ is the terminal object in the category of
$n$-retractible extensions of $\bfg$.  In particular, it is unique.
%item{(2)} If $n=1$, $\bfg$, $\bfg^\min(\emptyset)$ is just the
%direct product of groups $\bfg(\sigma)$ where $\sigma$ is a
%codimension 1 face of $S(\bfg)$.
\end{remark}

\begin{proof}  Let $\bfg^+_{dir}$ be the direct limit extension of
$\bfg$, as described in Proposition~\ref{18}. We recursively define
iterated unfoldings of $\bfg^+_{dir}$. For each simplex $\sigma_1$
of $|\bfg|$, define $d_{\sigma_1}\bfg^+_{dir}$ as previously.
Suppose that a complex of groups
$d_{\sigma_1,\ldots,\sigma_k}\bfg^+_{dir}$  has already been
defined.  For each simplex $\sigma_{k+1}$  of  the underlying
simplicial complex $|d_{\sigma_1,\ldots,\sigma_k}\bfg^+_{dir}|$, let
$d_{\sigma_1,\dots,\sigma_k,\sigma_{k+1}}\bfg^+_{dir}=
d_{\sigma_{k+1}}(d_{\sigma_1,\ldots,\sigma_k}\bfg^+_{dir})$.  Since
$\bfg$ is $n$-retractible, this allows us to define extended
complexes of groups $d_{\sigma_1,\ldots,\sigma_k}\bfg^+_{dir}$ for
any $k\leq n$. Due to Corollary~\ref{14}, all of these are extended
normal blocks of groups.

Note that each of the groups
$d_{\sigma_1,\ldots,\sigma_k}\bfg_{dir}(\emptyset)$, which we denote
$\tilG_{\sigma_1,\ldots,\sigma_k}$, is a subgroup of the direct
limit $\tilG$ of $\bfg$.  Let $N_k$ be the largest normal subgroup
of $\tilG$ which is contained in every subgroup
$\tilG_{\sigma_1,\ldots,\sigma_k}$, where $\sigma_1,\ldots,\sigma_k$
ranges over all {\it allowed} sequences of simplices (i.e., all sequences
for which $d_{\sigma_1,\ldots,\sigma_k}\bfg^+_{dir}$ was defined
above). If we set $\bfg_\min(\emptyset)=\tilG/N_n$, the resulting
extension is clearly $n$-retractible.  We need to show it is
minimal.

Let $\bfg^+_0$ be any $n$-retractible extension of $\bfg$. Clearly,
we have the canonical surjective homomorphism
$h:\tilG\to\bfg_0(\emptyset)$, which allows us to express
$\bfg_0(\emptyset)$ as the quotient $\tilG/\ker h$. To get the
homomorphism $\bfg_0(\emptyset)\to\bfg_\min(\emptyset)=\tilG/N_n$ as required
for minimality, we need to show that $\ker h<N_n$. By definition of
$N_n$, it is thus sufficient to show that $\ker
h<\tilG_{\sigma_1,\dots,\sigma_n}$ for all allowed sequences
$\sigma_1,\dots,\sigma_n$.

To prove the latter, we will show by induction on $k$ that $\ker
h<\tilG_{\sigma_1,\dots,\sigma_k}$, for all $1\le k\le n$. For
$k=1$, we have retractions $\tilG\to\bfg(\sigma_1)$ and
$\bfg_0(\emptyset)\to\bfg(\sigma_1)$ which commute with $h$, and
such that $\tilG_{\sigma_1}=\ker[\tilG\to\bfg(\sigma_1)]$ and
$d_{\sigma_1}\bfg_0(\emptyset)=\ker[\bfg_0(\emptyset)\to\bfg(\sigma_1)]$.
Consequently, we get a canonical homomorphism
$h_{\sigma_1}:\tilG_{\sigma_1}\to d_{\sigma_1}\bfg_0(\emptyset)$
with $\ker H_{\sigma_1}=\ker h$. In particular, $\ker
h<\tilG_{\sigma_1}$.

To proceed, observe that (due to Lemma~\ref{11} applied recursively)
the non-extended unfoldings $d_{\sigma_1,\dots,\sigma_k}\bfg_{dir}$
and $d_{\sigma_1,\dots,\sigma_k}\bfg_{0}$ coincide for all $1\le
k\le n$. Thus, for any
$\sigma_2\subset|d_{\sigma_1}\bfg_{dir}|=|d_{\sigma_1}\bfg_0|$ we
can repeat the above argument and get the homomorphism
$h_{\sigma_1,\sigma_2}:\tilG_{\sigma_1,\sigma_2}\to
d_{\sigma_1\sigma_2}\bfg_0(\emptyset)$ with $\ker
h_{\sigma_1,\sigma_2}=\ker h_{\sigma_1}=\ker h$. Repeating this
argument, we finally get the homomorphisms
$h_{\sigma_1,\dots,\sigma_n}:\tilG_{\sigma_1,\dots,\sigma_n}\to
d_{\sigma_1,\dots,\sigma_n} \bfg_0(\emptyset)$ with $\ker
h_{\sigma_1,\dots,\sigma_n}=\ker h$. It follows that $\ker
h<\tilG_{\sigma_1,\dots,\sigma_n}$. This proves statement~1.

To prove statement 2, recall that, by Definition~\ref{1}, the block
$|\bfg|$ is finite. If ${\bfg}(\sigma)$ is finite for each
$\sigma\neq \emptyset$, then it follows (by applying recursively
Lemmas~\ref{7} and~\ref{11}) that each underlying simplicial complex
$|d_{\sigma_1,\ldots,\sigma_k}\bfg_{dir}|$ is finite, and that for
every non-empty simplex
$\eta\subset|d_{\sigma_1,\ldots,\sigma_k}\bfg_{dir}|$ the group
$d_{\sigma_1,\ldots,\sigma_k}\bfg_{dir}(\eta)$ is finite. Hence each
$\tilG_{\sigma_1,\ldots,\sigma_k}$ has finite index in $\tilG$, and
there are finitely many such groups for each $k\leq n$.  If follows
that the intersection of all such subgroups is a subgroup of finite
index, and so each $N_k$ for $k\leq n$ is a finite index normal
subgroup of $\tilG$. This proves that $\bfg_\min=\tilG/N_n$ is a
finite group as claimed in statement~2.

Before proving statements 3~and~4, we first claim that each $N_k$ is
equal to the intersection of the groups
$\tilG_{\sigma_1,\ldots,\sigma_k}$.  To see this, it is useful to
change the indexing set.  For $\tau_1,\ldots,\tau_k$ a sequence of
simplices of $ d\bfg^+_{dir}$ of length at most $n$, define
$\sigma_1$ to be the image of $\tau_1$ in $|\bfg|=\tilG\backslash
d\bfg^+_{dir}$.  Assuming that $\sigma_1,\ldots,\sigma_{i-1}$ have
already been defined for some $i$ with $1<i\leq k$, define
$\sigma_i$ to be the image of $\tau_i$ in
$|d_{\sigma_1,\ldots,\sigma_{i-1}}\bfg^+_{dir}|=
\tilG_{\sigma_1,\ldots,\sigma_{i-1}}\backslash d\bfg^+_{dir}$.  Also
define $d_{\tau_1\ldots,\tau_k}\bfg^+_{dir}$ to be equal to
$d_{\sigma_1,\ldots, \sigma_k}\bfg^+_{dir}$, and define
$\tilG_{\tau_1\ldots,\tau_k}$ to be equal to
$\tilG_{\sigma_1,\ldots,\sigma_k}$.

If $x$ is a point of $ d\bfg^+_{dir}$ whose stabilizer is some
subgroup $H< \tilG$, and if $g$ is an element of $\tilG$, the
stabilizer of the point $g.x$ is equal to the conjugation
$gHg^{-1}$.  This observation and induction show that for each $g\in
\tilG$, for each $k\le n$ and for each sequence
$\tau_1,\ldots,\tau_k$ of simplices of $ d\bfg^+_{dir}$, we have
$$
\tilG_{g\tau_1,\ldots,g\tau_k}= g\cdot
\tilG_{\tau_1,\ldots,\tau_k}\cdot g^{-1}.
$$
Hence the intersection of the subgroups of the form
$\tilG_{\tau_1,\ldots,\tau_k}$, for fixed $k\leq n$, is a normal
subgroup of $\tilG$.  It follows that this intersection is equal to
$N_k$.

The above observation combined with the inclusions
$\tilG_{\sigma_1,\dots,\sigma_k}<\tilG_{\sigma_1,\dots,\sigma_k,\sigma_{k+1}}$
implies that $N_{k+1}<N_k$. To prove statement 3, we will show by
induction that for $0\le k<n$ the quotient groups $N_k/N_{k+1}$ and
$\tilG/N_{k+1}$ are $p$-groups of bounded exponent. Here we use
convention that $N_0=\tilG$.

For $k=0$, we need to show that $\tilG/N_1$ is a $p$-group of
bounded exponent. We know that $N_1$ is the intersection of the
groups of form $\tilG_{\sigma_1}$, which are the kernels of
retractions $\tilG\to\bfg(\sigma_1)$. Thus $\tilG/N_1$ embeds in the
product of the groups $\bfg(\sigma_1)$. Since the latter groups are
$p$-groups of bounded exponent, the assertion follows.

Now we suppose $\tilG/N_k$ is a $p$-group of bounded exponent and claim
that  the group  $N_k/N_{k+1}$ is too. To see that this is true, recall
that the groups $\tilG_{\sigma_1,\dots,\sigma_{k+1}}$ are the
kernels of the retractions $\tilG_{\sigma_1,\dots,\sigma_k}\to
d_{\sigma_1\dots\sigma_k}\bfg_{dir}(\sigma_{k+1})$. Note also that
$$
N_{k+1}=\bigcap\tilG_{\sigma_1,\dots,\sigma_{k+1}}=\bigcap[N_k\cap\tilG_{\sigma_1,\dots,\sigma_{k+1}}],
$$
and thus $N_{k+1}$ is equal to the intersection of the kernels of
the composed homomorphisms
$$
N_k\,\,\mapright{{\rm incl}}\,\,\tilG_{\sigma_1,\ldots,\sigma_k}
\,\,\mapright{r}\,\,
d_{\sigma_1,\ldots,\sigma_k}\bfg_{dir}(\sigma_{k+1}).
$$
Now, each of the groups
$d_{\sigma_1,\ldots,\sigma_k}\bfg_{dir}(\sigma_{k+1})$ canonically
embeds in the quotient $\tilG/N_k$ (because the latter is
a $k$-retractible extension of $\tilG$). Thus all these groups are
$p$-groups of finite exponent. As before, we see that $N_k/N_{k+1}$
embeds in the product of $p$-groups of bounded exponents, hence the
assertion. The fact that the quotient $\tilG/N_{k+1}$ is then also a
$p$-group of bounded exponent follows directly. This proves
statement 3.

The proof of statement 4 is similar to that of statement 3.
\end{proof}

\subsection{The construction and retra-products}\label{1.4}

Given any $n$,
we construct an $n$-retrac\-tible extended simplex of groups $\bfg^+$, over the
simplex $\Delta$ of arbitrary dimension, as follows:

{\bf 1.} We put the trivial group on the simplex $\Delta$. \pra

{\bf 2.} We put an arbitrary group $\bfg (s)$ on each codimension 1
face $s$ of $\Delta$. \pra

{\bf 3.} Suppose we  have already defined groups $\bfg(\eta)$ (and
inclusion maps between them) for faces of codimension strictly less
than $k$. Then for a face $\tau$ of codimension $k\le\dim\Delta$, the group
$\bfg(\tau)$ is the minimal $n$-retractible extension, as in Theorem
20, of the simplex of groups over the link simplex $\Delta_\tau$
made of the already defined groups $\bfg(\eta)$, via the canonical
correspondence between the faces of $\Delta_\tau$ and the faces
$\eta\subset\Delta$ containing $\tau$. \pra

{\bf 4.} Finally, we take as $\bfg(\emptyset)$ the minimal
$n$-retractible extension of the so far obtained simplex of groups
$\bfg$ over $\Delta$.

\begin{definition}[Retra-product]\label{21} We will call the group
$\bfg(\emptyset)$  of any simplex of groups $\bfg^+$ obtained as in
the construction above the \emph{$n$-retra-product} of the (finite)
family of groups $\bfg(s):s\in{\mathcal S}_\Delta$. Note that this
operation makes sense for any finite family of groups.
\end{definition}

Clearly, the groups $\bfg(\tau)$ obtained in the construction above
are all the $n$-retra-products of the corresponding families of
groups $\bfg(s):s\in{\mathcal S}_\Delta, s\supset\tau$.

Rephrasing Theorem~\ref{20} in the context of retra-products we get
the following properties of this operation, for an arbitrary natural
number $n$:
\begin{itemize}\item[(1)] the $n$-retra-product of  finite groups is finite;
\item[(2)] the $n$-retra-product of $p$-groups of bounded exponent is a $p$-group
of bounded exponent.
\end{itemize}

 It follows that the $n$-retra-product of finite $p$-groups
is a finite $p$-group. In particular, we get the following.

\begin{coroll}\label{22}    Let $\bfg$ be a non-extended simplex of
groups obtained as in the construction above, out of groups
$\bfg(s)$ being finite $p$-groups. Then all groups $\bfg(\tau)$ in
this simplex are finite $p$-groups.\end{coroll}

\begin{remark}  The construction of this subsection was first used
in~\cite{JS1}, in the $2$-retractible case.  For the purposes of this
paper, it suffices to consider the case when each codimension one
face of the simplex is assigned the cyclic group $\Z_p$ of order $p$. Even
though the construction of the $n$-retra-products is in principle
explicit, one rapidly loses track of the groups arising.

The $n$-retra-product of two groups $\Z_2$, occurring at faces of
codimension~2, is the dihedral group $D_{2^n}$ of order $2^{n+1}$.
We do not even know the orders of the $n$-retra-products of three copies
of $\Z_2$, except  the case  $n=2$ when this order is $2^{14}$.
\end{remark}

Note that if $\bfg$ is a non-extended simplex of groups
corresponding to $\bfg^+$, obtained by the construction above, then
its direct limit extension $\bfg^+_{dir}$ is different from
$\bfg^+$, and in particular the direct limit
$\bfg^+_{dir}(\emptyset)$ is different from the $n$-retra-product
$\bfg(\emptyset)$.  This motivates the following.

\begin{definition}[Free retra-product]\label{23} The direct limit of a
non-extended simplex of groups $\bfg$ obtained by the construction
above will be called the \emph{free $n$-retra-product} of the
(finite) family of the codimension~1 groups $\bfg(s):s\in{\mathcal
S}_\Delta$.
\end{definition}

In the next subsection we will deal with free $n$-retra-products of
finite groups, showing that for $n\ge2$ they are infinite, and for
$n\ge3$ they are non-elementary word-hyperbolic.

\subsection{Simplicial non-positive curvature, word-hyperbolicity,
and the proof of Theorem \ref{thm:monstrous}}\label{1.5}

To show that there are non-elementary word-hyperbolic groups arising
from the construction of the previous subsection, we will use
results of~\cite{JS2} concerning simplicial non-positive curvature.

Recall from~\cite{JS2} that the systole $sys(K)$ of a simplicial
complex $K$ is the smallest number of 1-simplices in any full
subcomplex of $K$ homeomorphic to the circle. A simplicial complex
$K$ is $k$-large if its systole and the systoles of links of all
simplices in $K$ are all at least $k$. Simplicial complexes whose
all links are $k$-large, for some fixed $k$, are the analogs of
metric spaces with curvature bounded above. If they are additionally
simply connected, we call them \emph{$k$-systolic complexes}. All
results of this subsection are corollaries to the following.

\begin{proposition}\label{24}  Suppose $\bfg^+$ is an $n$-retractible
extended simplex of groups. Then the development $d\bfg^+$ is
$2(n+1)$-large.
\end{proposition}

We skip the proof of the proposition until the end of the
subsection, first deriving (and making comments on) its
consequences. In particular, we show how this proposition,
together with the results of the previous subsection,
implies Theorem \ref{thm:monstrous}.

\begin{coroll}\label{25} Suppose $\bfg$ is a non-extended simplex of
groups whose links $\bfg_\sigma^+$ are $n$-retractible, and let
$\bfg^+_{dir}$ be the direct limit extension of $\bfg$.
\begin{itemize}
\item[(1)] If $n\ge2$ then the development $d\bfg^+_{dir}$ is contractible.
\item[(2)] If $n\ge2$ and the codimension~1 groups $\bfg(s)$ are non-trivial, then
$d\bfg^+_{dir}$ and the direct limit $\bfg_{dir}(\emptyset)$ are
both infinite.
\item[(3)] If $n\ge3$ then the 1-skeleton of the development $d\bfg^+_{dir}$,
equipped with the polygonal metric, is Gromov-hyperbolic.
\item[(4)] If $n\ge3$ and all the groups $\bfg(\sigma)$ are finite and nontrivial, then
the group $\bfg_{dir}(\emptyset)$ is non-elementary word-hyperbolic,
except when the underlying simplex is the 1-simplex and the vertex
groups are of order 2 (in which case it is the infinite dihedral group).
\item[(5)] If $n\ge 2$ and all the groups $\bfg(\sigma)$ are finite, then
the group $\bfg_{dir}(\emptyset)$ is virtually torsion-free.
\end{itemize}
\end{coroll}

\begin{proof} First, note that links in $d\bfg^+_{dir}$ are
isomorphic to the developments of the $n$-retractible link simplices
of groups $\bfg_\sigma$ (see Lemma 9(1)), and thus are
$2(n+1)$-large. Since $d\bfg^+_{dir}$ is the universal cover of
$\bfg$, it is simply connected and hence $2(n+1)$-systolic. For
$n\ge2$, this means that $d\bfg^+_{dir}$ is 6-systolic, and thus it
is contractible by~\cite[Theorem 4.1(1)]{JS2}. This proves (1).

To get assertion (2), note that it follows from~\cite[Proposition
18(2)]{JS2}, that the group $\bfg_{dir}(\emptyset)$ has a nontrivial
subgroup that acts on the development $d\bfg^+_{dir}$ freely. Thus
this development is a classifying space for this subgroup and, since
it has finite dimension, the subgroup has to be infinite. See
also~\cite[Corollary 4.3]{Os1},  for a more elementary argument.

For parts (3) and (4), note first that $n\ge3$ implies 8-systolicity
of the development $d\bfg^+_{dir}$. Part (3) follows then
from~\cite[Theorem 2.1]{JS2}.

Finally, under the assumptions of (4), the group
$\bfg_{dir}(\emptyset)$ acts on $d\bfg^+_{dir}$ properly
discontinuously and cocompactly. Thus it follows from (3) that
$\bfg_{dir}(\emptyset)$ is word-hyperbolic. If the underlying
simplex of $|\bfg|$ is 1-dimensional, the group acts geometrically
on the tree $d\bfg^+_{dir}$, and hence is virtually free non-abelian
(except the mentioned case). If the dimension of the underlying
simplex $|\bfg|$ is greater than $1$, the fact that the group is
non-elementary follows from~\cite[Theorem 5.6 and Remark 2 at the
end]{Os2}, where it is shown that this group has one end.

The non-elementariness above can be also shown directly, by noting
that the groups as in (4) contain as a subgroup the free product of
three nontrivial finite groups (namely codimension~1 groups in any
vertex unfolding of certain three pairwise disjoint sides). We do
not include the details of this argument.

In order to prove (5), recall that by Theorem \ref{20} there exists a finite extension
$\bfg_{min}(\emptyset)$ of $\bfg$. Since the canonical homomorphism
$\lambda: \bfg_{dir}(\emptyset) \to \bfg_{min}(\emptyset)$ is injective on the local groups,
the kernel $H$ of $\lambda$ acts freely on the development
$d\bfg^+_{dir}$ and the index $|\bfg_{dir}(\emptyset):H|$ is finite. Suppose that $H$ contains a non-trivial element $h$ of finite order $p$;
without loss of generality, we can assume that $p$ is prime.
The space $d\bfg^+_{dir}$ is finite-dimensional by definition, and
contractible by part (1), therefore, by Smith fixed point theory (see \cite[Thm. III.7.11]{bredon}), the
induced action of $\langle h \rangle$ on this space has a fixed point. This yields a contradiction with the freeness
of the action of $H$ on $d\bfg^+_{dir}$.
\end{proof}

Restricting to free $n$-retra-products, we immediately get the
following.

\begin{coroll}\label{26}  If $n\ge2$ then the free $n$-retra-product of
(at least two) nontrivial groups is infinite and virtually  torsion-free. If $n\ge3$ then the
free $n$-retra-product of (at least two) nontrivial finite groups is
non-elementary word-hyperbolic, except the product of two groups of
order 2.
\end{coroll}

\begin{remark}\label{27} Note that the construction in the present paper
gives new families of $k$-systolic groups, for arbitrary $k$,
different from those constructed in~\cite[Sections 17--20]{JS2}.
These are the free $n$-retra-products of arbitrary finite groups,
where $2(n+1)\ge k$.
\end{remark}

\begin{proof}[Proof of Theorem \ref{thm:monstrous}]

Let $\bfg$ be the 3-retractible $(n+1)$-simplex of groups obtained by
the construction of Subsection 3.4, with all codimension 1 groups
$\bfg(s)$ isomorphic to the cyclic group $\Z_p$ of order $p$.
Let $G$ be the direct limit of $\bfg$, i.e. the free 3-retra-product
of $n+2$ copies of $\Z_p$.
By Corollary 3.30, $G$ is then non-elementary word hyperbolic.

Choose a generator for each codimension 1 subgroup $\bfg(s)$ of $G$
(i.e. for each factor of the above free 3-retra-product).
Let $S$ be the set formed of these generators.
Then $S$ consists of $n+2$ elements, and it generates $G$
since the union $\cup\{ \bfg(s):s\in|\bfg| \}$ generates $G$.
For any proper subset $T\subset S$, the subgroup
$\langle T  \rangle<G$ generated by $T$ coincides with one
of the local groups of $\bfg$. More precisely,
for $\tau=\cap\{ s:s\in T \}$ we have $\langle T \rangle=\bfg(\tau)$.
By Corollary 3.25, $\langle T \rangle$ is then a finite $p$-group,
which completes the proof.
\end{proof}

It only remains now to prove Proposition \ref{24}. We will use the
following lemma in the proof.

\begin{lemma}\label{28}  Let $\gamma$ be a loop in the 1-skeleton of a
simplicial complex $K$ which has the minimum length $L$ amongst all
loops in the 1-skeleton of $K$ in the same free homotopy class. Then
any lift $\widetilde\gamma$ of $\gamma$ to the universal cover
$\widetilde K$ of $K$ has the property that it minimizes distance
measured in the 1-skeleton of $\widetilde K$ between any two points
whose distance in $\widetilde \gamma$ is at most $L$.
\end{lemma}

\begin{proof} Let $t$ be a deck transformation of $\widetilde K$ that
acts as a translation by the distance $L$ on the
subcomplex~$\widetilde \gamma$.  Suppose that~$\widetilde \gamma$
does not have the property claimed.  Then there exist vertices
$p$~and~$q$ of~$\widetilde\gamma$ such that $q$ lies on the segment
from $p$ to $tp$ and such that the distance between $p$ and $q$ in
the 1-skeleton of $\widetilde K$ is strictly smaller than the
distance between them in $\widetilde \gamma$.  Let $\alpha$ be the
segment of $\widetilde\gamma$ between $p$ and $q$, and let $\alpha'$
be a path of shorter length in the 1-skeleton of $\widetilde K$ from
$p$ to $q$. Denote also by $\beta$ be the segment of
$\widetilde\gamma$ between $q$ and $tp$.  Since $\widetilde K$ is
simply connected, $\alpha$ and $\alpha'$ are homotopic relative to
their endpoints. This homotopy, after projecting to $K$, yields a
homotopy between $\gamma$ and the loop obtained by projecting
$\alpha'\cup\beta$. Since the latter loop is strictly shorter, we
get a contradiction. \end{proof}

\begin{proof}[Proof of  Proposition~\ref{24}] The proof is by induction on $d$,
the dimension of the simplex $|\bfg|$, followed by induction on $n$.
For $d=0$ there is nothing to prove. Let $\bfg^+$ be an
$n$-retractible $d$-dimensional simplex of groups. Then links
$\bfg_\sigma^+$ are also $n$-retractible (see Remark 17(3)), and
thus by induction on $d$ their developments $d\bfg_\sigma^+$ are
$2(n+1)$-large. Thus the same holds for links of $d\bfg^+$.

By~\cite[Corollary 1.5]{JS2},  a simplicial complex $X$ with
$k$-large links is $k$-large iff the length of the shortest loop in
the 1-skeleton of $X$ which is homotopically nontrivial in $X$ is at
least $2(n+1)$. We thus need to show this for $X=d\bfg^+$.

Let $\gamma$ be a homotopically nontrivial polygonal loop in
$d\bfg^+$ of the shortest length. Now we start the induction on $n$.
It has been shown in~\cite[Proposition 4.3(3)]{JS1},  that the
development of any 1-retractible extended simplex of groups is a
flag complex. Hence the length of $\gamma$ is at least~4.  This
completes the case $n=1$ for all $d$.  Clearly an $n$-retractible
complex is $(n-1)$-retractible, and so by induction on $n$, the
length of $\gamma$ is at least $2n$.

Let $\widetilde X$ be the universal cover of $X=d\bfg^+$, let
$\widetilde\gamma$ be a lift of $\gamma$ to $\widetilde X$, and let
$v_0,\ldots,v_{n+1}$ be some consecutive vertices on
$\widetilde\gamma$.  By Lemma~\ref{28}, we see that
$\widetilde\gamma$ minimizes distances between the vertices
$v_i:0\le i\le n+1$ in the 1-skeleton of $\widetilde X$.  In
particular the distance between $v_0$ and $v_{n+1}$ is equal to
$n+1$.  Let $\delta$ denote the segment of $\widetilde\gamma$
between $v_0$ and $v_{n+1}$, and let $\delta'$ be the segment of
$\widetilde\gamma$ that starts at $v_{n+1}$ and projects to the
segment of $\gamma$ complementary to the projection of $\delta$.

By symmetry of $d\bfg^+$, we may assume that $v_1$ is a vertex of
the simplex $[|\bfg|,1]\subset d\bfg^+$. Consider the fundamental
domain $D\subseteq \widetilde X$ for the group
$\tilG_{v_1,\ldots,v_{n}}= d_{v_1,\ldots,v_{n}}(\emptyset)$ obtained
recursively as $d_{v_n}\dots d_{v_1}|\bfg|$ in the way described
just before Lemma 7. By Lemma 7, $D$ is a strict fundamental domain.
It contains the vertices $v_1,\dots,v_{n}$ in its interior (i.e.,
outside the boundary $\partial D$), and thus contains also $v_0$ and
$v_{n+1}$. We identify $D$ with the quotient
$\tilG_{v_1,\ldots,v_n}\backslash X$. We then look at the projection
of $\delta\cup\delta'$ to $D$. Since the distance in $\widetilde X$
between $v_0$ and $v_{n+1}$ is $n+1$, their distance in $D$ is also
$n+1$ (the distance in the quotient cannot increase, while that in
the subcomplex cannot increase). It follows that the length of the
projection of $\delta'$ to $D$ is at least $n+1$, and so the length
of $\gamma$ is at least $2(n+1)$.  \end{proof}

\subsection{The fixed point property}
We now pass to proving Theorem \ref{thm:helly}.
The proof will use the following lemma.
The proof of a related result, for contractible CW-complexes,
can be found in \cite{L}. We remind the reader that
``mod-$p$ acyclic'' means ``having the same mod-$p$ \v Cech cohomology as
a point''.

\def\bfg{{\bf G}}
\def\cals{{\cal S}}
\def\calsplus{{{\cal S}^+}}
\def\mapdown#1{\Big\downarrow\rlap{$\vcenter{\hbox{$\scriptstyle#1$}}$}}
\def\mapright#1{\smash{\mathop{\longrightarrow}\limits^{#1}}}
\def\min{{\rm min}}
\def\pra{{\par}}
\def\span{{\rm span}}

\begin{lemma}\label{lem:helly}  Fix an integer $n>0$,
let $Y_1,\ldots, Y_n$ be closed subspaces of a space $X$, let
$Y=\bigcup_{i=1}^n Y_i$ and let $A=\bigcap_{i=1}^nY_i$, the union and
intersection of the $Y_i$ respectively.  Suppose that for all subsets
$I\subset\{1,\dots,n\}$ with $1\leq |I|<n$, the intersection
$\bigcap_{i\in I}Y_i$ is mod-$p$ acyclic.  Then the reduced mod-$p$ \v
Cech cohomologies of $Y$ and $A$ are isomorphic, with a shift in
degree of $n-1$.  More precisely, for each $m$ there is an isomorphism
$\widetilde H^m(Y)\cong \widetilde H^{m-n+1}(A)$.
\end{lemma}

%\begin{lemma}\label{lem:helly}
%Let $Y_1,\ldots, Y_n$ be closed mod-$p$ acyclic
%subspaces of a space $X$ and suppose that for all proper subsets
%$I\subset\{1,\dots,n\}$
%the intersection $\bigcap_{i\in I}Y_i$   is also mod-$p$ acyclic.
%\pra
%\item{(i)}
%If $Y_1\cap\cdots\cap Y_n$ is mod-$p$ acyclic then so is $Y$.
%\pra
%\item{(ii)}
%If $Y_1\cap\cdots\cap Y_n$ is empty, then $Y$ has the same mod-$p$
%\v Cech cohomology as the $(n-2)$-sphere.
%
%\end{lemma}

\begin{proof}
In the case when $n=1$, we have that $Y=A$ and the assertion is
trivially true.  Now suppose that $n\geq 2$.  For $1\leq i\leq n-1$,
let $Z_i=Y_i\cap Y_n$, let $Z=\bigcup_i Z_i$, and let
$Y'=\bigcup_{i=1}^{n-1}Y_i$.  By definition, $Y-Y' = Y_n - Z$, and so
by the strong form of excision that holds for \v Cech cohomology (see
the end of section~3.3 of~\cite{hatcher}), it follows that
$H^*(Y,Y')\cong H^*(Y_n,Z)$.

By induction on $n$, we see that for each $m$, $\widetilde
H^m(Z)\cong \widetilde H^{m-n+2}(A)$.  Also by induction, we see that
$Y'$ is mod-$p$ acyclic, since $\bigcap_{i=1}^{n-1}Y_i$ is mod-$p$
acyclic by hypothesis.  Since $n\geq 2$, the hypotheses also imply
that $Y_n$ is mod-$p$ acyclic.  Hence the long exact sequence in
reduced cohomology for the pair $(Y_n,Z)$ collapses to isomorphisms, for
all~$i$, $\widetilde H^{i-1}(Z)\cong H^i(Y_n,Z)$.  Similarly, the
long exact sequence in reduced cohomology for the pair $(Y,Y')$
collapses to isomorphisms, for all $i$, $H^i(Y,Y')\cong\widetilde
H^i(Y)$.  Putting these isomorphisms together gives an isomorphism,
for all $i$, of reduced cohomology groups $\widetilde H^{i-1}(Z)\cong
\widetilde H^i(Y)$.  The claimed result follows.
\end{proof}

\begin{proof}[Proof of  Theorem~\ref{thm:helly}]
Suppose that $X$ is a $p$-acyclic
$G$-space of finite covering dimension.
Let $g_1,\ldots,g_{n+2}$ be the elements of $S$, and let $Y_i$ be the
points of $X$ that are fixed by $g_i$.
%By Smith theory, the
%fixed point set for the action of any finite $p$-group on $X$ is
%mod-$p$ acyclic (see Theorem~III.7.11 of \cite{Bre}).
By Smith theory, the fixed point set of the action of any finite
$p$-group on $X$ is mod-$p$ acyclic.  For an explicit
reference, see Theorem~III.7.11 of Bredon's book~\cite{bredon},
noting that ``finite covering dimension'' implies Bredon's hypothesis
``finitistic'' (see p.~133 of~\cite{bredon}), and that ``$X$ is mod-$p$
acyclic'' is equivalent to Bredon's hypothesis ``the pair
$(X,\emptyset)$ is a mod-$p$ \v Cech cohomology 0-disk''. Hence the
subspaces $Y_i$ satisfy the hypotheses of Lemma~\ref{lem:helly}.
The global fixed point set $A$ for the action of $G$ on $X$ is equal
to the intersection $A=\bigcap_{i=1}^{n+2}Y_i$.  If $A$ is not mod-$p$
acyclic, then for some $m\geq -1$, the reduced cohomology group
$\widetilde H^m(A)$ is non-zero.  By
Lemma~\ref{lem:helly}, it follows that the union
$Y=\bigcup_{i=1}^{n+2}Y_i$ has a non-vanishing reduced cohomology
group $\widetilde H^j(Y)$ for some $j\geq n$.  Hence the relative
cohomology group $H^{j+1}(X,Y)$ is non-zero for some $j\geq n$, and so
$X$ must have covering dimension at least $n+1$.
\end{proof}

%****************START HERE NEXT TIME.

%*******************************************************
\section{Finitely presented groups that fail to act}\label{sec:triv}
%******************************************************

In this section we establish the validity of Template
$\textbf{NA}_{fp}$.  We also show the triviality of actions on
manifolds of certain
groups, see  Lemma~\ref{lem:simpgp}. This immediately implies
Proposition~\ref{prop:simp} and provides the input to Template
$\textbf{NA}_{fp}$ that is needed to prove Theorem~\ref{thm:triv}.

\begin{definition}\label{def:rs} A sequence of groups and monomorphisms,
$(G_n;\xi_{n,j})\ (n\in\mathbb N, j=1,\dots,J)$, is called {\em{a
recursive system}} if
\begin{itemize}
    \item[(i)] each $G_n$ has a presentation $\langle \cA_n~\|~ R_n\rangle$
with $\cA_n$ finite and $\bigcup_n R_n\subset \cA^*$
recursively enumerable, where $\cA=\sqcup_n\cA_n$, and
    \item[(ii)] each
monomorphism $\xi_{n,j}:G_n\to G_{n+1}$ is defined by a set of words
$S_{n,j}=\{w_{n,j,a}\in \cA_{n+1}^*\mid a\in \cA_n\}$ such that
$w_{n,j,a}=\xi_{n,j}(a)$ in $G_{n+1}$, with
$\bigcup_{n,j}S_{n,j}\subset\cA^*$ recursively enumerable.
\end{itemize}

\end{definition}

We shall be interested only in sequences where, for each
sufficiently large integer $n$,
$G_{n+1}$ is generated by the union of the images of the
$\xi_{n,j}$. And in our applications we shall need only the case
$J=2$.

\begin{examples} \label{exy}
The following are recursive systems.

(1) Define $G_n=\SL(n,\Z)$, set $J=2$, and define
 $\xi_{n,1}$ and $\xi_{n,2}$ to be the embeddings
$\SL(n,\Z)\to\SL(n+1,\Z)$ defined by
$$\xi_{n,1}(M)=
\left( \begin{array}{cc} M & 0 \\ 0 & 1 \end{array}\right),~\mbox{ and }~
\xi_{n,2}(M)=\left( \begin{array}{cc} 1 & 0 \\ 0 & M \end{array}\right)
~\mbox{ for all } M \in \SL(n,\Z).$$

\medskip

(2) Writing $Alt(n)$ to denote the alternating group consisting of
even permutations of $\text{\bf n}=\{1,\dots,n\}$ and
taking $J=2$, define $G_n=Alt(n)$ and define
$\xi_{n,1},\xi_{n,2}:Alt(n)\to Alt(n+1)$ to be the embeddings induced
by the maps $I_{n,1}, I_{n,2}:\text{\bf{n}}\to\text{\bf{n+1}}$ defined
by $I_{n,1}:k\to k$ and $I_{n,2}:k\to k+1$.

\end{examples}

\begin{lemma}\label{trivlem} Let $(G_n;\xi_{n,j})\ (n\in\mathbb N,
j=1,\dots,J)$ be a recursive system of non-trivial groups and
monomorphisms $\xi_{n,j}:G_n\to G_{n+1}$, and
suppose that there
exists $n_0$ so that for each $n\geq n_0$,
$G_{n+1}$ is generated by $\bigcup_j\xi_{n,j}(G_n)$.  Then there
exists a finitely presented group $G_\omega$ which for each
$n\geq n_0$, contains two
isomorphic copies  of $G_n$ so that $G_\omega$ is the normal
closure of the union of these two subgroups.
\end{lemma}

%%We are now ready to construct our finitely presented example.
%Our result is the following
%\begin{thm} \label{thm:triv_act_man} There exists a finitely presented group $Q$ admitting non-trivial action on any finitely
%dimensional CAT(0)-manifold.
%\end{thm}

\begin{proof} By renumbering, we may assume that $n_0=1$.
Consider the free product \begin{equation} \label{eq:G}
G=\bigast_{n=1}^\infty G_n,
\end{equation} and for each $j\in J$ let $\theta_j:G\to G$ be the injective
endomorphism of $G$ whose restriction to $G_n$ is
$\xi_{n,j}$.

Now, let $H$ be the multiple HNN-extension of $G$ corresponding to
these endomorphisms:
$$H=\langle G, t_1,\dots,t_j~\|~t_jgt_j^{-1}=\theta_j(g),~g \in G, j\in
J \rangle~.$$
In the notation of Definition~\ref{def:rs}, this has presentation
$$
\langle \cA, t_1,\dots,t_j ~\|~ R_n \ (n\in\mathbb N),\
t_jat_j^{-1}w_{n,j,a}^{-1}\ (n\in
\mathbb N,\, j\in J, a\in \cA_n)\rangle.
$$
By hypothesis, this is a recursive presentation.
Moreover, since the images of the $\xi_{n,j}$ generate $G_{n+1}$, the group
$H$ is generated by the finite set $\cA_1\cup\{t_1,\dots,t_J\}$.

Now, by the Higman embedding theorem (see \cite[IV.7]{L-S}), $H$ can
be isomorphically embedded into a finitely presented group $B$.
Suppose that $B$ is generated by elements $x_1,\dots,x_l$. Without
loss of generality we can and do assume that each of $x_1,\dots,x_l$ has
infinite order. Indeed, to ensure this one can if necessary replace $B$ by the free
product $B*\Z$ of $B$ with the infinite cyclic group
generated by $z$, which is generated by the elements
$z,zx_1,\dots,zx_l$ of infinite order.  Choose an element of infinite
order $y \in G \le B$ and a subgroup $F \le G$ such that $F$ is free
of rank $l$ and
\begin{equation} \label{eq:F} F \cap \langle y \rangle=\{1\},~ F \cap
  \langle x_i \rangle =\{1\}~\mbox{
for}~i=1,\dots,l.
\end{equation} Such a choice is possible because $G$ is a free product
of infinitely many non-trivial groups and a cyclic
subgroup of $B$ can intersect at most one free factor non-trivially.

Consider, now, the iterated HNN-extension of $B$:
$$L=\langle B,s_1,\dots,s_l~\|~s_ix_is_i^{-1}=y,~i=1,\dots,l\rangle.$$

Let $\{f_1,\dots,f_l\}$ be free basis of $F$. By \eqref{eq:F} and
Britton's lemma (\cite[IV.2]{L-S}), the subgroup of $L$ generated by
$s_1,\dots,s_l$ and $F$ is freely generated by the elements
$s_1,\dots,s_l$, $f_1,\dots,f_l$. Let $L'$ be a copy of $L$ and let
$s'_1,\dots,s'_l,f'_1,\dots,f'_l$ denote the copies of the
corresponding elements. Finally we obtain the group
that we seek by defining
$$G_\omega=\langle L,L'~\|~s_i=f_i',f_i=s_i',~i=1,\dots,l\rangle.$$
The group $G_\omega$ is infinite and finitely presented by construction.

A key feature in our
construction is that
for all $k\le n$ the free factor $G_k$ of $G$ is conjugate in $H$
(and hence in $G_\omega$) to a subgroup of $G_n$, and for
$k>n$  the conjugates of $G_n$ by
positive words of length $k-n$ in the letters $t_j$  generate
$G_k$. Thus $G$ is the normal closure of each $G_n$. Likewise
$G'$ is the normal closure of $G_n'$.
All that remains is to
observe that $G_\omega$ is generated by the set
$\{x_1,\dots,x_l,s_1,\dots,s_l,x'_1,\dots,x'_l,s'_1,\dots,s'_l\}$,
each of whose elements is conjugate to an element of $G$
or $G'$.  Thus $G_\omega$ is the normal closure of $G_n\cup G_m'$
for every $n,m\ge 1$.
\end{proof}

The following theorem establishes the validity of Template $\text{\bf NA}_{fp}$.

\begin{thm}\label{trivAct} If the groups $G_n$ satisfy the conditions
of Template $\text{\bf NA}_{fp}$, then the finitely presented group
$G_\omega$ constructed in Lemma~\ref{trivlem} cannot act non-trivially
on any $X\in \mathcal{X}$.
\end{thm}

\begin{proof}
Suppose that $G_\omega$ acts on a space $X\in \mathcal{X}$.
Then $X\in \mathcal{X}_m$ for some $m\in\N$ and we have
a  homomorphism $\alpha: G_\omega \to {\rm{Homeo}}(X)$
that we want to prove is trivial. By hypothesis, there
is some $G_n$ that cannot
act non-trivially on $X$. Hence, in the notation of the preceding lemma,
$\alpha(G_n)=\alpha(G_n')=\{1\}$. Therefore the
kernel of $\alpha$ is the whole of $G_\omega$.
\end{proof}

\begin{lemma} \label{lem:simpgp}
Let $G$ be a simple group that contains a copy of $\Z_p^k$.
Then $G$ cannot act non-trivially on any mod-$p$ acyclic manifold $X$ of
dimension at most $n$, where $n=2k-1$ if $p$ is an odd prime and $n=k-1$ if
$p=2$.
\end{lemma}

%If $G$ is a simple group that contains a copy of
%$\Z_p^{n+1}$, then $G$ cannot act effectively by diffeomorphisms
%on any smooth $p$-acyclic manifold $X$ of dimension at most $n$.

\begin{proof}
Let $E$ be a subgroup of $G$ isomorphic to $\Z_p^k$.  It suffices to
show that $E$ cannot act effectively on $X$, since then the kernel of
the action of $G$ is not the trivial group and so must be equal to
$G$.

A proof that $E$ cannot act effectively on $X$ (with slightly weaker
hypotheses on $X$) is given by Bridson-Vogtmann in~\cite{Brid-Vogt},
who attribute the result to P.~A. Smith~\cite{pasmith}.  The result
can also be deduced from theorem~2.2 of~\cite{mannsu}, since a
standard result from Smith theory (Theorem~III.7.11 of~\cite{bredon})
implies that the $E$-fixed point set in $X$ is non-empty.

A more elementary proof that $E$ cannot act effectively can be
given if one restricts to the case when $X$ is smooth and $G$ acts by
diffeomorphisms.  First, fix a Riemannian metric on $X$ that is
compatible with the given smooth structure.  Next, average this metric
over the $E$-action to replace it by a metric for which $E$ acts on
$X$ by isometries.  As above, Smith theory ensures that there is a
point $x\in X$ fixed by $E$.  Taking derivatives at $x$ gives a linear
action of $E$ on the tangent space $T_x(X)$.  Now $E$ has no faithful
real representation of dimension at most $n$, and so some non-identity
element $h\in E$ acts trivially on $T_x(X)$.  Since the action of $E$
is by isometries, it follows that $h$ fixes an open ball around $x$.
But the fixed point set for any isometry is a closed submanifold of
$X$, and hence $h$ fixes all points of $X$.
\end{proof}

\noindent{\bf{Proof of Theorem \ref{thm:triv}}}
The preceding lemma tells us that, given any prime $p$
and positive integer $n$, any
alternating group $Alt(m)$ with $m$ sufficiently large cannot
act non-trivially on a $p$-acyclic
manifold of dimension less than $n$. It follows that the
group $G_\omega$ obtained by applying Theorem \ref{trivAct} to
the recursive system in Example \ref{exy}(2) cannot act
non-trivially on a $p$-acyclic manifold of
any dimension, for all primes $p$.

\medskip

\begin{remark}\label{rem:last}

For each fixed prime $p$, C. R\"over constructed a finitely presented
simple group containing, for each $n$, a copy of
$(\Z_p)^n$~\cite{Rov}.  By Lemma~\ref{lem:simpgp}, such a group cannot
act non-trivially on any mod-$p$ acyclic manifold.

Certain of the finitely presented simple groups introduced by Thompson
\cite{thomp} and Higman \cite{higman} contain a copy of every finite
group, so these cannot act on a mod-$p$ acyclic
manifold for any $p$.
\end{remark}

\end{document}